# TOEPLITZ OPERATORS AND
# WEIGHTED BERGMAN KERNELS

## Miroslav Engliš

ABSTRACT. For a smoothly bounded strictly pseudoconvex domain, we describe the boundary singularity of weighted Bergman kernels with respect to weights behaving like a power (possibly fractional) of a defining function, and, more generally, of the reproducing kernels of Sobolev spaces of holomorphic functions of any real order. This generalizes the classical result of Fefferman for the unweighted Bergman kernel. Finally, we also exhibit a holomorphic continuation of the kernels with respect to the Sobolev parameter to the entire complex plane. Our main tool are the generalized Toeplitz operators of Boutet de Monvel and Guillemin.

## 1. INTRODUCTION

Let $\Omega$ be a bounded strictly pseudoconvex domain in $\mathbf{C}^n$ with smooth boundary, and $r$ a defining function for $\Omega$. Thus $r$ is smooth on the closure $\overline{\Omega}$ of $\Omega$, $r < 0$ on $\Omega$, and $r = 0$, $\|\nabla r\| > 0$ on $\partial\Omega$. For simplicity of notation, we also use the "positively signed" defining function $\rho = -r$. It was then shown by Fefferman [12] that the Bergman kernel $K(x, y)$ of $\Omega$, when restricted to the diagonal $x = y$, has the form

$$(1) \qquad K(x,x) = \frac{a(x)}{\rho(x)^{n+1}} + b(x) \log \rho(x) \qquad \text{with } a, b \in C^\infty(\overline{\Omega}).$$

Further,

$$(2) \qquad a|_{\partial\Omega} = \frac{n!}{\pi^n} J[\rho],$$

where $J[\rho]$ is the Monge-Ampere determinant

$$(3) \qquad J[\rho] = (-1)^n \det \begin{bmatrix} \rho & \frac{\partial \rho}{\partial \rho} \\ \frac{\partial}{\partial \rho} & \frac{\partial}{\partial \overline{\partial} \rho} \end{bmatrix}$$

which is positive on $\partial\Omega$ in view of the strict pseudoconvexity of $\Omega$.

The formula (1) also extends to $K(x, y)$ for $x \neq y$, in the following way. We will say that a function $f \in C^\infty(\overline{\Omega} \times \overline{\Omega})$ is *almost-sesquianalytic* if $\partial f(x, y)/\partial \overline{x}$ as well as $\partial f(x, y)/\partial y$ vanish to infinite order on the diagonal $x = y$. For a real-valued

1991 *Mathematics Subject Classification.* Primary 32W25; Secondary 32A36, 32A25, 47B35.

Research supported by GA AV ČR grant no. A1019304 and Ministry of Education research plan no. MSM4781305904.

Typeset by $\mathcal{A}\mathcal{M}\mathcal{S}$-TEX

function $f(x)$ on $\Omega$, its *almost-sesquianalytic extension* is any almost-sesquianalytic function $f(x,y)$ on $\overline{\Omega} \times \overline{\Omega}$ such that $f(x,x) = f(x) \; \forall x \in \overline{\Omega}$ and

$$(4) \qquad \overline{f(x,y)} = f(y,x).$$

It is well known that such an extension always exists, and is unique up to functions vanishing on the diagonal to infinite order. Let $\rho(x,y)$ be an almost-sesquianalytic extension of $\rho$. From the Taylor expansion, one has

$$\rho(x,y) + \rho(y,x) - \rho(x) - \rho(y) = \sum_{j,k=1}^{n} \frac{\partial^2 r(x)}{\partial x_j \partial \overline{x}_k} (y_j - x_j)(\overline{y}_k - \overline{x}_k) + O(|x-y|^3).$$

By the strict pseudoconvexity of $\Omega$, it therefore follows that for $\epsilon > 0$ small enough

$$(5) \qquad 2 \operatorname{Re} \rho(x,y) \geq \rho(x) + \rho(y) + c|x-y|^2$$

if $x, y \in \overline{\Omega}$ and $|x-y| < \epsilon$. Adjusting $\rho(x,y)$ by a function which vanishes near the diagonal and satisfies the symmetry condition (4), we can achieve that (5) holds everywhere on $\overline{\Omega} \times \overline{\Omega}$. Then $\operatorname{Re} \rho(x,y) \geq 0$ for all $x, y \in \overline{\Omega}$, with equality occurring only if $x = y \in \partial \Omega$.

Throughout the rest of this paper, we will thus fix an almost-sesquianalytic extension $\rho(x,y)$ of $\rho = -r$ which satisfies (4) and (5) on $\overline{\Omega} \times \overline{\Omega}$.

In particular, $\log \rho(x,y)$ and $\rho(x,y)^\alpha$, for any $\alpha \in \mathbf{C}$, are thus well defined and smooth on $\Omega \times \Omega$.

The analogue of (1) for $x \neq y$ then asserts that there exist almost-sesquianalytic functions $a(x,y)$ and $b(x,y)$ on $\overline{\Omega} \times \overline{\Omega}$ such that

$$(6) \qquad K(x,y) = \frac{a(x,y)}{\rho(x,y)^{n+1}} + b(x,y) \log \rho(x,y)$$

for all $x, y \in \Omega$.

The aim of this paper is to establish a similar formula also for weighted Bergman kernels $K_w(x,y)$ — i.e. the reproducing kernels of the subspaces $L^2_{\mathrm{hol}}(\Omega, w)$ of all holomorphic functions in $L^2(\Omega, w)$ — for smooth positive weights $w$ on $\Omega$ which behave like a power of $\rho$, in the sense that

$$(7) \qquad w = \rho^\alpha e^g, \qquad \text{with } \alpha > -1 \text{ and } g \in C^\infty(\overline{\Omega}).$$

Furthermore, we also establish various generalizations of (6) for the reproducing kernels of Sobolev spaces of holomorphic functions (Sobolev-Bergman kernels).

If $\alpha = m$ is a positive integer and $\log \frac{1}{w}$ is strictly plurisubharmonic, it is possible to derive the weighted analogue of (6) by an argument due to Forelli and Rudin [13] and Ligocka [25]. Namely, if $\widetilde{\Omega}$ denotes the Hartogs domain

$$\widetilde{\Omega} = \{(x,t) \in \Omega \times \mathbf{C}^m : \|t\|^{2m} < w(x)\}$$

in $\mathbf{C}^{n+m}$, then $K_w(x,y)$ is the restriction of the (unweighted) Bergman kernel $\widetilde{K}((x,t),(y,s))$ of $\widetilde{\Omega}$ to the hyperplane $s = t = 0$:

$$(8) \qquad K_w(x,y) = \frac{\pi^m}{m!}\, \widetilde{K}((x,0),(y,0)).$$

Moreover, the hypothesis that $\alpha$ be a positive integer ensures that $\widetilde{\Omega}$ is smoothly bounded, while the strict plurisubharmonicity of $\log\frac{1}{w}$ implies that it is strictly pseudoconvex. Finally, the function $r((x,t)) = \|t\|^2 - w(x)^{1/m}$ is clearly a defining function for $\widetilde{\Omega}$. Applying Fefferman's expansion (6) to $\widetilde{\Omega}$, it thus follows from (8) that

$$(9) \qquad K_w(x,y) = \frac{a(x,y)}{\rho(x,y)^{n+m+1}} + b(x,y)\,\log\rho(x,y) \qquad \forall x,y \in \Omega.$$

Also,

$$a(x,x) = \frac{(n+m)!}{m!\pi^n}\, J[\rho](x) \qquad \forall x \in \partial\Omega.$$

Using the localization lemma of Fefferman ([12], Lemma 1 on page 6), it can be shown that (9) in fact remains in force even if $\log\frac{1}{w}$ is assumed to be just plurisubharmonic (not necessarily strictly); however, the argument breaks down if $\log\frac{1}{w}$ is not plurisubharmonic, or if $\alpha$ is not an integer, or if $\alpha = 0$ and $g \not\equiv 0$.

Our first main result is the following.

**Theorem A.** *For $w$ of the form* (7),

$$(10) \qquad K_w(x,y) = \begin{cases} \dfrac{a(x,y)}{\rho(x,y)^{n+\alpha+1}} + b(x,y)\,\log\rho(x,y) & \text{if } \alpha \in \mathbf{Z}, \\[3mm] \dfrac{a(x,y)}{\rho(x,y)^{n+\alpha+1}} + b(x,y) & \text{if } \alpha \notin \mathbf{Z}, \end{cases}$$

*for all $x,y \in \Omega$ and some almost-sesquianalytic $a,b \in C^\infty(\overline{\Omega} \times \overline{\Omega})$. Moreover,*

$$(11) \qquad a(x,x) = \frac{\Gamma(n+\alpha+1)}{\Gamma(\alpha+1)\,\pi^n}\, \frac{J[\rho](x)}{e^{g(x)}} \qquad \text{for } x \in \partial\Omega.$$

It should be noted that the leading term (11) was given by Hörmander [20] for $\alpha = 0$, and by Ligocka [25] for $g \equiv 0$ and any $\alpha > -1$. The fact that $K_w(x,y)$ is $C^\infty$ on $\overline{\Omega} \times \overline{\Omega}$ minus the boundary diagonal was shown for $\alpha \neq 0$ by Peloso [29]. (For $\alpha = 0$ it can be proved by a straightforward modification of Kerzman's original argument [22] for the unweighted Bergman kernel.)

Our second main result concerns, in a sense, an extension of (6) and (10) to $\alpha \leq -1$: namely, to the reproducing kernels of *Sobolev spaces* of holomorphic functions.

Recall that the ordinary Sobolev space $W^s(\Omega)$, on a smoothly bounded domain $\Omega \subset \mathbf{R}^n$, with $s$ a nonnegative integer, consists of all functions on $\Omega$ whose distributional derivatives of orders up to $s$ belong to $L^2$, with the norm

$$(12) \qquad \|f\|_s := \left( \sum_{|\nu| \leq s} \|D^\nu f\|_{L^2(\Omega)}^2 \right)^{1/2}$$

where the summation is over all multiindices $\nu$ of length not exceeding $s$, and $D^\nu$ stands for the corresponding differentiation. Equivalently, $W^s(\Omega)$ is the completion of $C^\infty(\overline{\Omega})$ with respect to the norm (12). The closure, in $W^s(\Omega)$, of the subset $\mathcal{D}(\Omega)$ of all functions in $C^\infty(\overline{\Omega})$ whose support is a compact subset of $\Omega$ is denoted by $W_0^s(\Omega)$. For negative integers $s$, $W^s(\Omega)$ is defined as the space of all distributions supported on $\overline{\Omega}$ which are bounded as linear functionals on $W_0^s(\Omega)$ — that is, loosely speaking, as the dual of $W_0^s(\Omega)$ with respect to the $L^2$-pairing

$$\langle f, g\rangle_{L^2(\Omega)} := \int_\Omega f\,\overline{g}$$

(the integral being taken with respect to the Lebesgue measure). Finally, for non-integer $s$ the spaces $W^s(\Omega)$ are defined by interpolation.

Alternatively, one can first define the Sobolev spaces on the whole $\mathbf{R}^n$, for any real order $s$, as the spaces of tempered distributions $f$ whose Fourier transform satisfies

$$(13) \qquad \|f\|_s' := \left(\int_{\mathbf{R}^n} (1+\|\xi\|^2)^s |\hat{f}(\xi)|^2\right)^{1/2} < \infty.$$

One the defines $W^s(\Omega)$ as the subspace of all $f \in W^s(\mathbf{R}^n)$ which, in some appropriate sense, are "supported" on $\Omega$. On the level of norms, this amounts to taking instead of (12) the (equivalent) norms

$$(14) \qquad \|f\|_s' := \left(\sum_{|\nu|\le s} \frac{s!}{\nu_1!\dots\nu_n!} \|\partial^\nu f\|_{L^2(\Omega)}^2\right)^{1/2}.$$

See e.g. [27] for more details on Sobolev spaces.

In general, functions in the Sobolev spaces of fractional order are rather difficult to characterize explicitly; fortunately, the situation is much better if $\Omega$ is a bounded strictly pseudoconvex domain in $\mathbf{C}^n$ with smooth boundary and we restrict attention to holomorphic functions — that is, to subspaces $W_{\mathrm{hol}}^s$ of holomorphic functions in $W^s$.

First of all, it is known ([26], Remark 1 on p.31) that for $s < \frac{1}{2}$, the space $W_{\mathrm{hol}}^s(\Omega)$ coincides with $L_{\mathrm{hol}}^2(\Omega, |r|^{-2s})$, with equivalent norms. Thus our Theorem A can be interpreted as a result about reproducing kernels of "Sobolev-Bergman spaces" (cf. [17]); however, unfortunately, since spaces with equivalent norms may have kernels with wildly different boundary singularities (see §8.1 below for an example), this does not tell us anything about the singularity of the kernel of $W_{\mathrm{hol}}^s(\Omega)$ with the norm (12) or (14).

Secondly, it is a result of Beatrous [1] that for any $s \in \mathbf{R}$ and $m$ a nonnegative integer, a holomorphic function $f$ belongs to $W_{\mathrm{hol}}^s$ if and only if $\partial^\nu f \in W_{\mathrm{hol}}^{s-m}$ for all multiindices $\nu$ with $|\nu| \le m$; further, $\|f\|_s$ is equivalent to the norm

$$\left(\sum_{|\nu|\le m} \|\partial^\nu f\|_{W^{s-m}}^2\right)^{1/2}.$$

In particular, combining this with the facts from the previous paragraph, we see that if $m > s - \frac{1}{2}$ then $f \in W^s_{\mathrm{hol}}$ if and only if $\partial^\nu f \in L^2_{\mathrm{hol}}(\Omega, |r|^{2m-2s})$ whenever $|\nu| \leq m$, and $\|f\|_s$ is equivalent to the norm

$$(15) \qquad \left( \sum_{|\nu| \leq m} \|\partial^\nu f\|^2_{L^2(\Omega, |r|^{2m-2s})} \right)^{1/2}.$$

(We will actually obtain a new proof of this equivalence in Sections 6 and 7 below as a byproduct.) Again, this equivalence does not tell us anything about the boundary behaviour of the kernels, except for the leading term.

Finally, it is also a result of [1] that one need not consider all derivatives $\partial^\nu$ in (15), but just "complex normal" derivatives: namely, let $\mathcal{D}$ be the holomorphic vector field on $\Omega$ "orthogonal" to $\partial r$, i.e. defined by $\mathcal{D}f = \langle df, \partial r \rangle$. Then a holomorphic function $f$ belongs to $W^s_{\mathrm{hol}}$ if and only if $\mathcal{D}^j f \in W^{s-m}$ for all $0 \leq j \leq m$, and $\|f\|_s$ is equivalent to

$$\left( \sum_{j=0}^m \|\mathcal{D}^j f\|^2_{W^{s-m}} \right)^{1/2}.$$

Again, for $m > s - \frac{1}{2}$ one can further replace here $W^{s-m}(\Omega)$ by $L^2_{\mathrm{hol}}(\Omega, |r|^{2m-2s})$:

$$(16) \qquad \left( \sum_{j=0}^m \|\mathcal{D}^j f\|^2_{L^2(\Omega, |r|^{2m-2s})} \right)^{1/2}.$$

(We will also obtain an independent proof of this equivalence as a byproduct in Section 7.)

Our second main result is the following. We actually expect it to be true for any $m > s - \frac{1}{2}$, but our proofs seem to work only for $m > 2s - 1$.

**Theorem B.** *Let $K^{(s)}(x, y)$ be the reproducing kernel of the holomorphic Sobolev space $W^s_{\mathrm{hol}}(\Omega)$ with respect to the norm corresponding to (13) (see Section 7 for the precise definition), or with respect to the norm (15) or (16) (for some nonnegative integer $m > 2s - 1$). Then*

$$(17)$$

$$K^{(s)}(x, y) = \begin{cases} \dfrac{a(x, y)}{\rho(x, y)^{n+1-2s}} + b(x, y) & \text{if } n + 1 - 2s \notin \mathbf{Z}, \\[2ex] \dfrac{a(x, y)}{\rho(x, y)^{n+1-2s}} + b(x, y) \log \rho(x, y) & \text{if } n + 1 - 2s \in \mathbf{Z}_{>0}, \\[2ex] \dfrac{a(x, y)}{\rho(x, y)^{n+1-2s}} \log \rho(x, y) + b(x, y) & \text{if } n + 1 - 2s \in \mathbf{Z}_{\leq 0}, \end{cases}$$

*on $\Omega \times \Omega$ for some almost-sesquianalytic functions $a, b \in C^\infty(\overline{\Omega} \times \overline{\Omega})$. Further, the values of $a(x, x)$ on the boundary are given by Corollary 21, Theorem 8 and Theorem 9, respectively.*

Notice that, in particular, for $s = \frac{1}{2}$ we recover the Szegö kernel, as the reproducing kernel of the Hardy space $H^2(\partial\Omega) \simeq W^{1/2}_{\mathrm{hol}}(\Omega)$.

Our third, and final, main result concerns holomorphic dependence of $K^{(s)}$ on the parameter $s$.

**Theorem C.** *Let $K^{(s)}(x, y)$ be the reproducing kernel of $W^s_{\text{hol}}(\Omega)$ with respect to the norm corresponding to* (13) *(see again Section 7 for the precise definition). Then the $K^{(s)}(x, y)$ extends to a holomorphic function of $x, \overline{y}, s$ on $\Omega \times \Omega \times \mathbb{C}$.*

To describe the idea of the proofs, recall that for any function $\phi \in L^\infty(\Omega)$, the Toeplitz operator $\mathbf{T}_\phi$ with symbol $\phi$ is the (bounded linear) operator on $L^2_{\text{hol}}(\Omega)$ defined by

$$(18) \qquad \mathbf{T}_\phi f := \mathbf{\Pi}(\phi f),$$

where $\mathbf{\Pi} : L^2(\Omega) \to L^2_{\text{hol}}(\Omega)$ is the orthogonal projection (the Bergman projection). It is immediate that if $\phi > 0$ on $\Omega$, then $\mathbf{T}_\phi$ is positive definite, hence, in particular, injective; thus there exists a (densely-defined, unbounded self-adjoint) inverse $\mathbf{T}_\phi^{-1}$. For $x \in \Omega$, let us write $K_x := K(\cdot, x)$ for the Bergman kernel at the point $x$, and similarly let $K_{w,x} := K_w(\cdot, x)$. The main ingredient in our proof is the simple observation that for any positive continuous weight function $w \in L^\infty(\Omega)$, $\mathbf{T}_w$ extends to a bounded operator from $L^2_{\text{hol}}(\Omega, w)$ into $L^2_{\text{hol}}(\Omega)$, and the kernels $K_x$ then belong to the range of $\mathbf{T}_w$ and

$$(19) \qquad K_{w,x} = \mathbf{T}_w^{-1} K_x.$$

One can also define Toeplitz operators even for some unbounded symbols $\phi$, as (unbounded) densely defined operators on $L^2_{\text{hol}}(\Omega)$: for instance, if $\phi \in L^2$ then (18) makes sense for all $f \in H^\infty(\Omega)$; and

$$(20) \qquad \mathbf{T}_\phi f(x) = \langle \phi f, K_x \rangle = \int_\Omega \phi\, f\, \overline{K_x}$$

makes sense for all $f \in H^\infty$ even if $\phi$ is only in $L^1$, since $K_x \in L^\infty(\Omega)$. In all these cases, at least for $w$ of the form (7), the formula (19) still remains in force.

Ideas very reminiscent of (19) have appeared in the paper of Bell [3].

There are also Toeplitz operators on the Hardy space $H^2(\partial\Omega)$, that is, the subspace in $L^2(\partial\Omega)$ of functions which are boundary values of holomorphic functions in $\Omega$ (alternatively, the closure in $L^2(\partial\Omega)$ of $C^\infty_{\text{hol}}(\overline{\Omega})$). Namely, for any $q \in L^\infty(\partial\Omega)$ the Toeplitz operator on $H^2(\partial\Omega)$ with symbol $q$ is defined by

$$(21) \qquad T_q f := \Pi(qf)$$

where $\Pi : L^2(\partial\Omega) \to H^2(\partial\Omega)$ is the orthogonal projection (the Szegö projection). More generally, following Boutet de Monvel [6], one can define a "generalized" Toeplitz operator $T_Q$ by the formula

$$(22) \qquad T_Q f := \Pi(Qf)$$

for any pseudodifferential operator $Q$ on $\partial\Omega$. The previous operators (21) are the special case when $Q$ is the operator of multiplication by the function $q$.

It is now a remarkable fact that for any $\phi \in C^\infty(\overline{\Omega})$ and any $Q$, both $\mathbf{T}_\phi$ and $T_Q$ map $C^\infty_{\text{hol}}(\overline{\Omega})$ into itself, and in fact for a given $\phi$ there exists a $Q$ such that, up to a "negligible" error term, $\mathbf{T}_\phi = T_Q$ on $C^\infty_{\mathrm{hol}}(\overline{\Omega})$. Further, if $\phi$ vanishes on $\partial\Omega$ to order $m$, then $Q$ can be chosen to be of order $-m$ (see [14], Theorem 9.1). We show that all this remains in force also for $\phi = w$ with $w$ of the form (7); the pseudo-differential operator $Q$ is then of (possibly fractional) order $-\alpha$. Combining this with the formula (19), the result of Theorem A follows by Boutet de Monvel's and Guillemin's calculus of the "generalized" Toeplitz operators (cf. their book [7], especially the Appendix) and the microlocal description of the singularity (6) as given by Boutet de Monvel and Sjöstrand [8]. For Theorem B, there is a formula analogous to (19) but with the role of $T_w$ taken over by a suitable (pseudo)differential operator; it turns out that this operator again coincides (on $C^\infty_{\mathrm{hol}}(\overline{\Omega})$) with a generalized Toeplitz operator $T_Q$ for some $Q$, and the result follows in a similar way as before. Finally, for Theorem C one combines the above facts with the existence, proved by Seeley [33], of the complex powers of an arbitrary positive elliptic pseudodifferential operator.

The formula (19), as well as its analogue needed for Theorem B, are established in Section 2. Section 3 overviews the necessary material from [6] and [7] on generalized Toeplitz operators. Theorem A is proved in Section 4. The part of Theorem B concerning the norms (15) and (16) is proved in Section 5; some further variations on this theme (including, as a byproduct, a new proof of the equivalence of (15) with the ordinary Sobolev norm, as well as extensions of this equivalence to spaces of harmonic functions) occur in Section 6. The part of Theorem B concerning the norm (13) and Theorem C are established in Section 7, together with a new proof and generalization to harmonic functions of the equivalence of the norms (16). The final Section 8 is occupied by miscellaneous concluding comments and remarks.

A different proof of Theorem A was obtained by G. Komatsu (personal communication), using Fefferman's original method from [12]; in fact, he was even able to handle weights $w$ with certain logarithmic-type singularities at the boundary. The present author hopes very much to see his work published soon.

The author also thanks Professor Louis Boutet de Monvel for illuminating him on several issues concerning generalized Toeplitz operators.

## 2. TOEPLITZ OPERATORS AND WEIGHTED BERGMAN KERNELS

Let $w$ be any positive, continuous and integrable weight function on $\Omega$, and let $\Lambda$ denote the operator $w^{1/2}\Pi$ on $L^2(\Omega)$. The domain of $\Lambda$ is $(L^2(\Omega) \ominus L^2_{\mathrm{hol}}(\Omega)) \oplus (L^2_{\mathrm{hol}}(\Omega, w) \cap L^2_{\mathrm{hol}}(\Omega))$; since $w$ is assumed to be integrable, the second summand contains the space $C^\infty_{\mathrm{hol}}(\overline{\Omega})$ of all functions in $C^\infty(\overline{\Omega})$ holomorphic on $\Omega$, hence, in particular, $\Lambda$ is densely defined. Being essentially the restriction to the closed subspace $L^2_{\mathrm{hol}}(\Omega)$ of the closed operator $f \mapsto w^{1/2}f$ of multiplication by $w^{1/2}$ on $L^2(\Omega)$, $\Lambda$ is also closed. Thus $\Lambda^*\Lambda$ is self-adjoint. Since $\Pi$ is bounded, $(\Pi w^{1/2})^* = w^{1/2}\Pi = \Lambda$, so $\Lambda^* = (\Pi w^{1/2})^{**} = \overline{\Pi w^{1/2}}$, where the last bar denotes closure (see e.g. [31], Chapter VIII, §1); thus

$$\Lambda^*\Lambda = \overline{\Pi w^{1/2}}\, w^{1/2}\Pi \supset \Pi w \Pi.$$

We declare this to be, *by definition*, the Toeplitz operator $T_w$. It is clear that for $w \in L^\infty(\Omega)$, this coincides with the operator defined by (18). Also, for $w \in L^2(\Omega)$, the domain of $T_w$ contains $H^\infty(\Omega)$, the space of bounded holomorphic functions on $\Omega$.

For any $f \in \operatorname{dom} T_w$, we have

$$T_w f(x) = \Lambda^* \Lambda f(x) = \langle \Lambda^* \Lambda f, K_x \rangle = \langle \Lambda f, \Lambda K_x \rangle$$
$$= \langle w^{1/2} f, w^{1/2} K_x \rangle = \int_\Omega w f \overline{K_x},$$

since $K_x \in C^\infty_{\mathrm{hol}}(\overline{\Omega}) \subset \operatorname{dom} \Lambda$ by Kerzman's theorem. Thus we see that, indeed, (20) holds for any $f \in \operatorname{dom} T_w$ and any $x \in \Omega$.

Since $w$ is continuous and positive, it is well known that the reproducing kernel $K_w(x,y)$ of $L^2_{\mathrm{hol}}(\Omega, w)$ exists and

$$(23) \qquad\qquad |f(x)| \le \|f\|_w \, \|K_{w,x}\|_w \qquad \forall f \in L^2_{\mathrm{hol}}(\Omega, w),$$

where $\| \cdot \|_w$ stands for the norm in $L^2(\Omega, w)$.

**Proposition 1.** *Assume that $K_x \in \operatorname{Ran} T_w$. Then (19) holds, i.e.*

$$K_{w,x} = T_w^{-1} K_x.$$

*Proof.* Let $K_x = T_w h$. Then for any $f \in L^2_{\mathrm{hol}}(\Omega) \cap L^2_{\mathrm{hol}}(\Omega, w) = \operatorname{dom} \Lambda$,

$$\langle f, K_{w,x} \rangle_w = f(x) = \langle f, \Lambda^* \Lambda h \rangle = \langle \Lambda f, \Lambda h \rangle$$
$$= \langle w^{1/2} f, w^{1/2} h \rangle = \int_\Omega w f \overline{h} = \langle f, h \rangle_w.$$

Since the set of all such $f$ is dense in $L^2_{\mathrm{hol}}(\Omega, w)$ (e.g. because it contains $C^\infty_{\mathrm{hol}}(\overline{\Omega})$), it follows that $h = K_{w,x}$.  $\square$

This simple observation can be generalized, as follows. Let, quite generally, $\Lambda$ be a closed, densely defined and injective operator on $L^2_{\mathrm{hol}}(\Omega)$. Introduce a new norm and inner product on $\operatorname{dom} \Lambda$ by

$$(24) \qquad\qquad \|f\|_\Lambda^2 := \|\Lambda f\|^2, \qquad \langle f, g \rangle_\Lambda := \langle \Lambda f, \Lambda g \rangle.$$

By von Neumann's theorem, the operator $\Lambda^* \Lambda$ is self-adjoint, and its square root $T = (\Lambda^* \Lambda)^{1/2}$ satisfies $\operatorname{dom} T = \operatorname{dom} \Lambda$, $\operatorname{Ran} T = \operatorname{Ran} \Lambda^*$ and $\|Tf\| = \|\Lambda f\|$ (see e.g. [10], Section XII.7). Thus we can write (24) equivalently as

$$(25) \qquad\qquad \|f\|_\Lambda^2 = \|Tf\|^2, \qquad \langle f, g \rangle_\Lambda = \langle Tf, Tg \rangle.$$

Let $\mathcal{H}_\Lambda$ denote the completion of $\operatorname{dom} \Lambda = \operatorname{dom} T$ with respect to this norm. (If $\Lambda^{-1}$ is bounded, then $\mathcal{H}_\Lambda$ coincides with $\operatorname{dom} \Lambda$; if $\Lambda$ is bounded, then $L^2_{\mathrm{hol}}(\Omega) \hookrightarrow \mathcal{H}_\Lambda$ continuously.)

**Proposition 2.** (i) *The evaluation functional* $f \mapsto f(x)$ *is continuous on* $\operatorname{dom}\Lambda$ *with respect to the* $\Lambda$-*norm, i.e.*

$$|f(x)| \leq c_x \, \|\Lambda f\| \qquad \forall f \in \operatorname{dom}\Lambda, \tag{26}$$

*if and only if*

$$K_x \in \operatorname{Ran}\Lambda^* (= \operatorname{Ran}T). \tag{27}$$

(ii) *If* (27) *is satisfied for all* $x \in \Omega$, *then the elements of* $\mathcal{H}_\Lambda$ *can be identified with functions on* $\Omega$ *and* $\mathcal{H}_\Lambda$ *admits the reproducing kernel*

$$K_\Lambda(x,y) \equiv K_{\Lambda,y}(x) = \langle T^{-1}K_y, T^{-1}K_x \rangle. \tag{28}$$

(iii) *If* $K_x$ *belongs not only to* $\operatorname{Ran}\Lambda^* = \operatorname{Ran}T$ *but even to* $\operatorname{Ran}\Lambda^*\Lambda = \operatorname{Ran}T^2$, *then*

$$K_{\Lambda,x} = (\Lambda^*\Lambda)^{-1}K_x. \tag{29}$$

*Proof.* (i) If $K_x = Th$, then

$$|f(x)| = |\langle f, Th \rangle| = |\langle Tf, h \rangle| \leq \|Tf\| \, \|h\| = \|\Lambda f\| \, \|h\|$$

for all $f \in \operatorname{dom}\Lambda = \operatorname{dom}T$, so (26) is satisfied with $c_x = \|h\|$. Conversely, if (26) holds, then

$$\Lambda f \mapsto f(x) = \langle f, K_x \rangle$$

is a well-defined bounded linear functional on $\operatorname{Ran}\Lambda$; extending it to the whole space by continuity, it follows that there exists $h \in L^2_{\text{hol}}(\Omega)$ such that

$$\langle f, K_x \rangle = \langle \Lambda f, h \rangle \qquad \forall f \in \operatorname{dom}\Lambda.$$

This means that $h \in \operatorname{dom}\Lambda^*$ and $\Lambda^* h = K_x$. So $K_x \in \operatorname{Ran}\Lambda^* = \operatorname{Ran}T$.

(iii) If $\{f_n\}$ is a sequence in $\operatorname{dom}\Lambda$ Cauchy in the $\Lambda$-norm, (26) implies that $\{f_n(x)\}$ is a Cauchy sequence in $\mathbb{C}$; thus the evaluation functionals $f \mapsto f(x)$ extend continuously from $\operatorname{dom}\Lambda$ to all of $\mathcal{H}_\Lambda$. Further, this extension still satisfies (26). It follows that elements of $\mathcal{H}_\Lambda$ can be viewed as functions on $\Omega$, and that $\mathcal{H}_\Lambda$ admits a reproducing kernel $K_\Lambda(x,y)$. If $\{f_n\}, \{g_n\}$ are sequences in $\operatorname{dom}\Lambda$ converging in $\mathcal{H}_\Lambda$ to $K_{\Lambda,x}$ and $K_{\Lambda,y}$ respectively, then by (25) $\{Tf_n\}$ is a Cauchy sequence in $L^2_{\text{hol}}(\Omega)$, whence $Tf_n \to F$ for some $F \in L^2_{\text{hol}}(\Omega)$. Since $\langle f, K_x \rangle = f(x) = \langle f, K_{\Lambda,x} \rangle_\Lambda = \lim \langle f, f_n \rangle_\Lambda = \lim \langle Tf, Tf_n \rangle = \langle Tf, F \rangle$ for any $f \in \operatorname{dom}T$, we have $K_x = T^*F$, or $F = T^{-1}K_x$. Similarly $Tg_n \to G = T^{-1}K_y$. Thus

$$K_\Lambda(x,y) = \langle K_{\Lambda,x}, K_{\Lambda,y} \rangle_\Lambda = \lim \langle g_n, f_n \rangle_\Lambda = \lim \langle Tg_n, Tf_n \rangle$$
$$= \langle G, F \rangle = \langle T^{-1}K_y, T^{-1}K_x \rangle,$$

proving (28).

(iii) This is immediate from (28), and also is easy to check directly: if $K_x = \Lambda^*\Lambda h$, then for any $f \in \operatorname{dom}\Lambda$

$$\langle f, h \rangle_\Lambda = \langle \Lambda f, \Lambda h \rangle = \langle f, \Lambda^*\Lambda h \rangle = f(x) = \langle f, K_{\Lambda,x} \rangle_\Lambda.$$

Since $\operatorname{dom}\Lambda$ is dense in $\mathcal{H}_\Lambda$, it follows that $h = K_{\Lambda,x}$. $\square$

If $K_x$ belongs to $\operatorname{Ran}\Lambda^*$ but not to $\operatorname{Ran}\Lambda^*\Lambda$, then $K_{\Lambda,x}$ does not belong to $\operatorname{dom}\Lambda$ but only to the part of $\mathcal{H}_\Lambda$ obtained by completion; in that case the formula (29) still holds, in the following sense.

**Proposition 3.** *The operators $\Lambda$ and $T : \operatorname{dom}\Lambda \to L^2_{\mathrm{hol}}(\Omega)$ extend by continuity to unitary isomorphisms, still denoted $\Lambda$ and $T$, from $\mathcal{H}_\Lambda$ onto $L^2_{\mathrm{hol}}(\Omega)$, and if (27) holds then $T^2 K_{\Lambda,x} = K_x$ (or $K_{\Lambda,x} = T^{-2} K_x$).*

*Proof.* The first part is immediate from (24) and (25). If (27) holds, then we have seen in the proof of part (ii) of Proposition 2 that for any sequence $f_n$ in $\operatorname{dom}\Lambda$ converging to $K_{\Lambda,x}$ in the $\mathcal{H}_\Lambda$-norm, $T f_n \to T^{-1} K_x \in L^2_{\mathrm{hol}}(\Omega)$ in $L^2_{\mathrm{hol}}(\Omega)$; by the first part of the current proposition, $T^{-1} K_x = Th$ for some $h \in \mathcal{H}_\Lambda$ and $f_n \to h$ in $\mathcal{H}_\Lambda$. Since $f_n \to K_{\Lambda,x}$, it follows that $h = K_{\Lambda,x}$ and $K_x = T^2 h = T^2 K_{\Lambda,x}$.  $\square$

We conclude this section by proving a result which will be needed further on. For $\alpha$ an integer this is well-known, and easily proved since $w$, being then $C^\infty$ on $\overline{\Omega}$, is a multiplier of $W^s(\Omega)$ and the Bergman projection $\mathbf{\Pi}$ is bounded from $W^s$ into $W^s_{\mathrm{hol}}$, for any $s \in \mathbf{R}$. For $\alpha$ noninteger, a proof seems to be needed.

**Proposition 4.** *Let $w$ be of the form (7). Then $T_w$ maps $C^\infty_{\mathrm{hol}}(\overline{\Omega})$ into itself.*

*Proof.* We use an idea of Bell [2]. Let $\phi_j$, $j = 1, \ldots, n$, be functions in $C^\infty(\overline{\Omega})$ such that $\phi_j \geq 0$ on $\overline{\Omega}$, $\sum_j \phi_j = 1$ near $\partial\Omega$, and $\partial_j r \neq 0$ on the support of $\phi_j$. Here $\partial_j r(z)$ is an abbreviation for $\frac{\partial r(z)}{\partial z_j}$ and $r(z)$ is the defining function for $\Omega$ from the Introduction.

First of all, we claim that for any $h \in L^2_{\mathrm{hol}}(\Omega) \cap C(\overline{\Omega})$ and $v \in L^\infty(\Omega)$,

$$(30) \qquad \int_\Omega \overline{h}\, \partial_j(|r|^{\alpha+1} v) = 0.$$

To see this, fix $\epsilon > 0$ and integrate by parts to obtain

$$\int_{r < -\epsilon} \overline{h}\, \partial_j(|r|^{\alpha+1} v) = \int_{r = -\epsilon} \overline{h}\, |r|^{\alpha+1}\, v\, n_j - \int_{r < -\epsilon} |r|^{\alpha+1} v\, \partial_j \overline{h}.$$

(Here $n_j$ is the appropriate component of the outward unit normal to $\partial\Omega$.) The second integral vanishes since $\overline{\partial}_j h = 0$, while the first integral is bounded by

$$\operatorname{area}(\{r = -\epsilon\})\, \cdot\, \|h\|_\infty\, \cdot\, \|v\|_\infty\, \cdot\, \epsilon^{\alpha+1}$$

which tends to zero as $\epsilon \searrow 0$ since $\alpha > -1$.

Let us now apply this observation to $h = K_x$ and $v = e^g f \phi_j / \partial_j r$, where $f \in C^\infty(\overline{\Omega})$. We get

$$\int_\Omega \partial_j\Big(\phi_j \frac{|r|}{\partial_j r}\, f w\Big) \overline{K_x} = 0.$$

Consequently,

$$\mathbf{\Pi} w f(x) = \int_\Omega w\, f\, \overline{K_x}$$
$$= \int_\Omega \Big[ wf - \sum_{j=1}^n \partial_j\Big(\phi_j \frac{r}{\partial_j r}\, wf\Big) \Big] \overline{K_x}.$$

By the Leibniz rule, the expression in the square brackets equals

$$(1 - \sum_j \phi_j)wf - rwf \sum_j \partial_j\left(\frac{\phi_j}{\partial_j r}\right) - rf \sum_j \frac{\phi_j}{\partial_j r}\partial_j w - rw \sum_j \frac{\phi_j}{\partial_j r}\partial_j f.$$

Since $\partial_j w = \partial_j(|r|^\alpha e^g) = |r|^\alpha e^g \partial_j g - \alpha|r|^{\alpha-1}e^g \partial_j r = (\partial_j g + \alpha\frac{\partial_j r}{r})w$, this equals

$$(1 - \sum_j \phi_j)wf - rw\left[f \sum_j \partial_j\left(\frac{\phi_j}{\partial_j r}\right) + f \sum_j \frac{\phi_j}{\partial_j r}\partial_j g + \sum_j \frac{\phi_j}{\partial_j r}\partial_j f\right] - \alpha wf$$
$$= |r|wf_1 - \alpha wf,$$

where

$$f_1 = \frac{1 - \sum_j \phi_j}{|r|}f + \sum_j \left(\partial_j\left(\frac{\phi_j}{\partial_j r}\right) + \frac{\phi_j}{\partial_j r}\partial_j g\right)f + \sum_j \frac{\phi_j}{\partial_j r}\partial_j f \in C^\infty(\overline{\Omega}).$$

Thus

$$(\alpha+1)\mathbf{\Pi}wf(x) = \int_\Omega |r|wf_1\overline{K_x} = \mathbf{\Pi}|r|wf_1(x);$$

that is,

$$(\alpha+1)\mathbf{\Pi}wf = \mathbf{\Pi}|r|wf_1.$$

Since $|r|w$ is a bounded function, it follows first of all that $\mathbf{\Pi}wf \in L^2_{\mathrm{hol}}(\Omega)$, that is, $f \in \mathrm{dom}\, T_w$. Further, repeating the above process with $|r|w$ in the place of $w$, we obtain successively functions $f_m \in C^\infty(\overline{\Omega})$, $m = 1, 2, \ldots$, such that

$$\frac{\Gamma(\alpha+m+1)}{\Gamma(\alpha+1)} \mathbf{\Pi}wf = \mathbf{\Pi}|r|^m wf_m.$$

Since $|r|^m wf_m = |r|^{m+\alpha+1}e^g f_m$ belongs to $W^s(\Omega)$ for $s < m + \alpha + \frac{3}{2}$ and $\mathbf{\Pi}$ maps $W^s(\Omega)$ into $W^s_{\mathrm{hol}}(\Omega)$ for any real $s$, it follows that $\mathbf{\Pi}wf \in W^s_{\mathrm{hol}}(\Omega)$ for all $s < m + \alpha + \frac{3}{2}$. As $m$ is arbitrary, $\mathbf{\Pi}wf \in \cap_{s\in\mathbf{R}}W^s_{\mathrm{hol}}(\Omega) = C^\infty_{\mathrm{hol}}(\overline{\Omega})$, completing the proof. $\square$

## 3. Generalized Toeplitz operators

Denote by $\eta$ the restriction to $\partial\Omega$ of the 1-form $\mathrm{Im}\,\partial r = (\partial r - \overline{\partial}r)/2i$. The strict pseudoconvexity of $\Omega$ is reflected in the fact that $\eta$ is a contact form, i.e. $\eta \wedge (d\eta)^{n-1}$ determines a nonvanishing volume element on $\partial\Omega$, or, equivalently, the half-line bundle

$$\Sigma := \{(x,\xi) \in T^*(\partial\Omega) : \xi = t\eta_x,\ t > 0\}$$

is a symplectic submanifold of $T^*(\partial\Omega)$. Equip $\partial\Omega$ with a measure with smooth positive density, and let $L^2(\partial\Omega)$ be the Lebesgue space with respect to this measure, and $L^2_{\mathrm{hol}}(\partial\Omega) = H^2(\partial\Omega)$ the subspace of nontangential boundary values of functions holomorphic in $\Omega$. We will also denote by $W^s(\partial\Omega)$, $s \in \mathbf{R}$, the Sobolev spaces on $\partial\Omega$, and by $W^s_{\mathrm{hol}}(\partial\Omega)$ the corresponding subspaces of nontangential boundary values of functions holomorphic in $\Omega$. (Thus $W^0(\partial\Omega) = L^2(\partial\Omega)$ and $W^0_{\mathrm{hol}}(\partial\Omega) = H^2(\partial\Omega)$.)

Unless otherwise specified, by a pseudodifferential operator or Fourier integral operator (ΨDO or FIO for short) on $\partial\Omega$ we will always mean an operator which is "regular" or "classical", i.e. in any local coordinate system the total symbol has an asymptotic expansion

$$p(x,\xi) \sim \sum_{j=0}^{\infty} p_{m-j}(x,\xi),$$

where $p_{m-j}$ is $C^\infty$ for $\xi \neq 0$, and positively homogeneous of degree $m-j$ with respect to $\xi$; here $j$ runs through nonnegative integers, but $m$ can be any complex number, and the symbol "$\sim$" means that the difference between $p$ and $\sum_{j=0}^{k-1} p_{m-j}$ should belong to the Hörmander class $S_{0,1}^{\operatorname{Re} m-k}$, for each $k=0,1,2,\dots$. If $P,Q$ are ΨDOs we write $P \sim Q$ if $P-Q$ is smoothing (i.e. of degree $-\infty$, or, equivalently, given by a $C^\infty$ Schwartz kernel).

For $P$ a ΨDO of order $m$ on the compact manifold $\partial\Omega$, the *generalized Toeplitz operator* $T_P : W_{\mathrm{hol}}^m(\partial\Omega) \to H^2(\partial\Omega)$ is defined as

$$T_P = \Pi P \Pi,$$

where $\Pi : L^2(\partial\Omega) \to H^2(\partial\Omega)$ is the orthogonal projection (the Szegö projector). Actually, $T_P$ maps continuously $W^s(\partial\Omega)$ into $W_{\mathrm{hol}}^{s-m}(\partial\Omega)$, for any $s \in \mathbf{R}$, because $\Pi$ is bounded on $W^s(\partial\Omega)$ for any $s \in \mathbf{R}$ (see [8]).

Microlocally, generalized Toeplitz operators have the following structure. Let $(x,q)$ denote the variable in $\mathbf{R}^n \times \mathbf{R}^{n-1} \simeq \mathbf{R}^{2n-1}$, and let $(\xi,v)$ be the dual variable. We identify $T^*\mathbf{R}^n$ with the symplectic cone $\Sigma_0 \subset T^*\mathbf{R}^{2n-1}$ defined by $y = v = 0$, and set

$$D_j = \frac{\partial}{\partial v_j} + v_j |D_x|, \qquad j = 1, \dots, n-1.$$

Let $H_0 : C^\infty(\mathbf{R}^n) \to C^\infty(\mathbf{R}^{2n-1})$ be the Hermite operator

$$H_0\phi(x,y) = (2\pi)^{-n} \int_{\mathbf{R}^n} e^{ix\cdot\xi - \frac{1}{2}\|\xi\|\,y\cdot y} \left(\frac{\|\xi\|}{2\pi}\right)^{(n-1)/4} \hat\phi(\xi)\,d\xi$$

where we write $x \cdot \xi$ for $\sum_j x_j \xi_j$, and the hat denotes Fourier transform. Then it follows from [5] and [8] that $\Pi$ admits the following microlocal description: for any $z_0 \in \partial\Omega$, there exists a canonical map $\Phi$ from a conic open set $U \subset T^*\mathbf{R}^{2n-1} \setminus \{0\}$ to a conic neighbourhood $V$ of $(z_0, \eta_{z_0}) \in \Sigma \subset T^*\Omega \setminus \{0\}$, whose restriction defines a symplectic isomorphism $\chi : \Sigma \cap U \to \Sigma \cap V$. There exists an elliptic FIO $F$, defined in $V$ modulo smoothing operators, associated with $\chi$, which transforms the left ideal of ΨDOs generated by the $D_j$ into the left ideal generated by the components of $\overline{\partial}_b$. Set $A \sim H_0^* F^* F H_0$ (this is an elliptic positive ΨDO) and $H \sim F H_0 A^{-1/2}$ (this is a FIO with complex phase, cf. [28]). Then $H^*H \sim I$, $HH^* \sim \Pi$, and for any ΨDO $Q$ on $\partial\Omega$,

$$(31) \qquad T_Q = \Pi Q \Pi \sim HPH^* \text{ near } z_0, \text{ with } \quad P \sim H^*QH \sim H^*T_Q H.$$

In fact the map $T_Q \mapsto P \sim H^*T_Q H$ is onto. It follows as a corollary that the generalized Toeplitz operators have the following properties.

(P1) Generalized Toeplitz operators form an algebra which is, modulo smoothing operators, locally isomorphic to the algebra of ΨDOs on $\mathbf{R}^n$.

(P2) In fact, for any $T_P$ there exists a $\Psi$DO $Q$ such that $T_P = T_Q$ and $Q\Pi = \Pi Q$.

(P3) It can happen that $T_P = T_Q$ for two different $\Psi$DOs $P$ and $Q$. If $\mathrm{ord}(P) > \mathrm{ord}(Q)$, then the restriction of the principal symbol $\sigma(P)$ of $P$ to $\Sigma$ vanishes. If $\mathrm{ord}(P) = \mathrm{ord}(Q)$, then the restrictions of the principal symbols $\sigma(P)$ and $\sigma(Q)$ to the cone $\Sigma$ coincide. Thus we can define unambiguously the order of $T_Q$ as $\min\{\mathrm{ord}(P) : T_P = T_Q\}$, and the symbol of $T_Q$ as $\sigma(T_Q) := \sigma(Q)|_\Sigma$ if $\mathrm{ord}(Q) = \mathrm{ord}(T_Q)$.

(P4) The order and the symbol are multiplicative: $\mathrm{ord}(T_Q T_{Q'}) = \mathrm{ord}(T_Q) + \mathrm{ord}(T_{Q'})$ and $\sigma(T_Q T_{Q'}) = \sigma(T_Q)\sigma(T_{Q'})$.

(P5) If $\mathrm{ord}(T_P) \leq 0$, then $T_P$ is a bounded operator on $L^2(\partial\Omega)$; if $\mathrm{ord}(T_P) < 0$, then it is even compact.

(P6) If $\mathrm{ord}(T_P) = \mathrm{ord}(T_Q) = k$ and $\sigma(T_P) = \sigma(T_Q)$, then $\mathrm{ord}(T_P - T_Q) \leq k-1$. In particular, if $T_P \sim T_Q$, then there exists a $\Psi$DO $Q' \sim Q$ such that $T_P = T_{Q'}$.

(P7) We will say that a generalized Toeplitz operator $T_P$ of order $m$ is elliptic if $\sigma(T_P)$ does not vanish. Then $T_P$ has a parametrix, i.e. there exists a Toeplitz operator $T_Q$ of order $-m$, with $\sigma(T_Q) = \sigma(T_P)^{-1}$, such that $T_P T_Q \sim I \sim T_Q T_P$.

We refer to the book [7], especially the Appendix, and to the paper [6] (which we have loosely followed in this section) for the proofs and additional information on generalized Toeplitz operators.

## 4. WEIGHTED BERGMAN KERNELS

Let now $\mathbf{T}_w$, with the weight function $w$ as in (7), be our Toeplitz operator on $L^2_{\mathrm{hol}}(\Omega)$ considered in Section 2. By Proposition 4 $\mathbf{T}_w$ maps $C^\infty_{\mathrm{hol}}(\overline{\Omega})$ into itself. In particular, since $C^\infty_{\mathrm{hol}}(\overline{\Omega})$ is contained and dense also in $H^2(\partial\Omega)$, we can view $\mathbf{T}_w$ as a (densely defined) operator on $H^2(\partial\Omega)$. Our first goal is to identify this operator with a certain generalized Toeplitz operator from the preceding section. For $\alpha$ an integer (so that $w \in C^\infty(\overline{\Omega})$), this was done by Guillemin ([14], Theorem 9.1), using ideas from [6]; we will show that essentially the same argument works also here.[1]

Following [6], let, quite generally, $\Lambda$ be an elliptic positive $\Psi$DO on $\partial\Omega$. Similarly as in Section 2, we define the Hilbert space $W^\Lambda$ as the completion of $C^\infty(\partial\Omega)$ with respect to the norm

$$(32) \qquad \|f\|^2_\Lambda := \langle \Lambda f, f \rangle = \int_{\partial\Omega} \overline{f}\, \Lambda f.$$

(We will usually use the notation $\langle \cdot, \cdot \rangle$ for the inner products both in $L^2(\Omega)$ and $L^2(\partial\Omega)$ — it is clear from the context which of the two is meant.) We further denote by $W^\Lambda_{\mathrm{hol}}$ the subspace of boundary values of functions holomorphic in $\Omega$, and by $\Pi_\Lambda : W^\Lambda \to W^\Lambda_{\mathrm{hol}}$ the orthogonal projector. Finally, for any $\Psi$DO $Q$ on $\partial\Omega$, we can again define a "generalized Toeplitz operator" associated to $\Lambda$ by

$$T^{(\Lambda)}_Q := \Pi_\Lambda Q \Pi.$$

---

[1]For $\Omega$ the unit disc or, equivalently, the upper half-plane in $\mathbf{C}$, where everything is much more tractable, Guillemin's proof was worked out more explicitly by Peng and Wong [30].

As shown in Section 1.d of [6], these operators possess pretty much the same microlocal description as the ordinary generalized Toeplitz operators $T_Q$ in the preceding section. Namely, with the notations $A, F, H_0$ and $H$ introduced before (31), setting

$$A_\Lambda \sim H_0^* F^* \Lambda F H_0 \sim A^{1/2} H^* \Lambda H A^{1/2},$$
$$H_\Lambda \sim F H_0 A_\Lambda^{-1/2} \sim H A^{1/2} A_\Lambda^{-1/2},$$

we have $H^* \Lambda \Lambda H_\Lambda \sim I$, $H_\Lambda H^* \Lambda \Lambda \sim \Pi_\Lambda$ (so that, in particular, modulo smoothing operators $H_\Lambda$ is an isomorphism from $L^2(\mathbf{R}^n)$ onto $W_{\text{hol}}^\Lambda$), and

$$
\begin{aligned}
(33) \qquad T_Q^{(\Lambda)} &= \Pi_\Lambda Q \Pi \sim H P' H^* \sim T_{P'}, \\
&\text{where } P' \sim H^* \Pi_\Lambda Q H \sim A^{1/2} A_\Lambda^{-1} A^{1/2} H^* \Lambda Q H.
\end{aligned}
$$

Let $\mathbf{K}$ denote the "Poisson extension operator" solving the boundary value problem

$$(34) \qquad \overline{\partial}^* \overline{\partial} \mathbf{K} u = 0 \quad \text{on } \Omega, \qquad \mathbf{K} u|_{\partial\Omega} = u.$$

(Thus $\mathbf{K}$ acts from functions on $\partial\Omega$ into functions on $\Omega$. Note that $\overline{\partial}^* \overline{\partial}$ is essentially the ordinary Laplace operator, only the boundary conditions are different; in particular, the kernel of $\overline{\partial}^* \overline{\partial}$ contains all holomorphic functions.) It is not difficult to see that $\mathbf{K}$ is actually continuous from $L^2(\partial\Omega)$ into $L^2(\Omega)$, and its range coincides with the subspace $L^2_{\text{harm}}(\Omega)$ of all harmonic functions. The adjoint $\mathbf{K}^*$ is thus continuous from $L^2(\Omega)$ to $L^2(\partial\Omega)$. The composition $\mathbf{K}^* \mathbf{K}$ is known to be an elliptic positive $\Psi$DO $\Lambda_0$ on $\partial\Omega$ of order $-1$, and we have

$$(35) \qquad \Lambda_0^{-1} \mathbf{K}^* \mathbf{K} = I_{\partial\Omega},$$

while

$$(36) \qquad \mathbf{K} \Lambda_0^{-1} \mathbf{K}^* = \mathbf{\Pi}_{\text{harm}},$$

the orthogonal projection in $L^2(\Omega)$ onto $L^2_{\text{harm}}(\Omega)$. Comparing (36) with (34), we also see that the restriction of $\Lambda_0^{-1} \mathbf{K}^*$ to $L^2_{\text{harm}}(\Omega)$ is the operator $\gamma$ of "taking the boundary values" of a harmonic function.

Thus, the operator $\mathbf{T}_w = \mathbf{\Pi} w|_{L^2_{\text{hol}}(\Omega)}$, when viewed as an operator on $H^2(\partial\Omega)$, is simply $\gamma \mathbf{\Pi} w \mathbf{K}|_{H^2(\partial\Omega)}$. Now by (35) and (36),

$$
\begin{aligned}
(37) \qquad \gamma \mathbf{\Pi} w \mathbf{K} &= \gamma \mathbf{\Pi} \mathbf{\Pi}_{\text{harm}} w \mathbf{K} \\
&= \gamma \mathbf{\Pi} \mathbf{K} \Lambda_0^{-1} \mathbf{K}^* w \mathbf{K}.
\end{aligned}
$$

Observe now that for $\Lambda = \Lambda_0$, (32) becomes

$$\|f\|_{\Lambda_0}^2 = \langle \mathbf{K}^* \mathbf{K} f, f \rangle_{L^2(\partial\Omega)} = \|\mathbf{K} f\|_{L^2(\Omega)}^2.$$

It follows that $u \mapsto \mathbf{K}u$ is an isometry, with inverse $\gamma$, identifying the space $W^{\Lambda_0}$ with $L^2_{\mathrm{hol}}(\Omega)$. In particular,

$$\gamma \mathbf{\Pi K} = \Pi_{\Lambda_0}.$$

Thus we can continue (37) by

$$\gamma \mathbf{\Pi} w \mathbf{K} = \Pi_{\Lambda_0} \Lambda_0^{-1} \Lambda_w,$$

where

$$\Lambda_w := \mathbf{K}^* w \mathbf{K}.$$

It was shown by Boutet de Monvel [4] that $\Lambda_w$, for $w$ as in (7), is an elliptic positive $\Psi$DO on $\partial\Omega$ of order $-\alpha - 1$ and with symbol

$$\sigma(\Lambda_w)(x, \xi) = \frac{\Gamma(\alpha + 1)}{2\|\xi\|^{\alpha+1}} \, e^{g(x)} \, \|\partial r(x)\|^\alpha.$$

(In fact, [4] only discussed integer $\alpha > -1$, but the case of noninteger $\alpha$ can be treated in the same manner, see the computation on the bottom of p. 256 and the remarks on the top of p. 257 in [6].) Consequently, by (33),

$$\gamma \mathbf{\Pi} w \mathbf{K}|_{H^2(\partial\Omega)} = T^{(\Lambda_0)}_{\Lambda_0^{-1} \Lambda_w} \sim T_Q,$$

with $Q = A^{1/2} A_{\Lambda_0}^{-1} A^{1/2} H^* \Lambda_w H$ a $\Psi$DO having the same order and symbol on $\Sigma$ as $\Lambda_0^{-1} \Lambda_w$.

We have thus arrived at the following proposition.

**Proposition 5.** *Let $w$ be a weight of the form (7). Viewed as on operator on $C^\infty_{\mathrm{hol}}(\overline{\Omega}) \subset H^2(\partial\Omega)$, the Toeplitz operator $\mathbf{T}_w$ on $L^2_{\mathrm{hol}}(\Omega)$ then coincides, modulo a smoothing operator, with the generalized Toeplitz operator $T_Q$ for some $\Psi$DO $Q$ on $\partial\Omega$, $\mathrm{ord}(Q) = -\alpha$, $\sigma(Q)|_\Sigma = \Gamma(\alpha + 1) \|\xi\|^{-\alpha} \|\partial r\|^\alpha e^g$.*

We are now ready to prove our first main result.

*Proof of Theorem A.* Let $Q$ be the operator from the last proposition. Since $\sigma(T_Q) > 0$, $T_Q$ is elliptic, so by (P7) there exists a $\Psi$DO $P$ (a parametrix) such that $\mathrm{ord}(P) = -\mathrm{ord}(Q) = \alpha$,

$$(38) \qquad \sigma(T_P) = \sigma(T_Q)^{-1} = \frac{\|\xi\|^\alpha e^{-g}}{\Gamma(\alpha + 1)\|\partial r\|^\alpha},$$

and $T_P \sim T_Q^{-1}$. By the property (P2) from Section 3, we can also assume without loss of generality that $P$ commutes with $\Pi$: $\Pi P = P\Pi$, so that $T_P$ is just the restriction of $P$ to $H^2(\partial\Omega)$.

Since $\mathbf{T}_w$ (or, more precisely, $\gamma \mathbf{T}_w \mathbf{K}$; we will drop $\gamma$ and $\mathbf{K}$ for the rest of this section) $\sim T_Q$, we thus have $\mathbf{T}_w^{-1} \sim T_P$, or $\mathbf{T}_w^{-1} \sim$ (the restriction of $P$ to $H^2(\partial\Omega)$). Using now our formula (19), we thus see that

$$K_{w,y} \sim T_P K_y = P K_y$$

(where "$\sim$" means that the two sides differ by a function in $C^\infty(\overline{\Omega})$). Now according to the main result of [8], there exists a classical symbol

$$b(x,y,t) \sim \sum_{j=0}^{\infty} t^{n-j} \, b_j(x,y), \qquad b_j \in C^\infty(\overline{\Omega} \times \overline{\Omega}),$$

such that the ordinary (unweighted) Bergman kernel satisfies

$$(39) \qquad K(x,y) = \int_0^\infty e^{-t\rho(x,y)} \, b(x,y,t) \, dt.$$

Thus

$$(40) \qquad PK_y \sim \int_0^\infty \sum_{j=0}^{\infty} t^{n-j} \, P\big[e^{-t\rho(\cdot,y)} \, b_j(\cdot,y)\big] \, dt.$$

By the standard symbol calculus for $\Psi$DOs (see, for instance, Theorem 4.2 in Hörmander [19]), we have quite generally (i.e. for any $\Psi$DO $P$ of order $\alpha$; here $P$ applies to the $x$ variable)

$$t^{n-j} P\big[e^{-t\rho(x,y)} \, b_j(x,y)\big] = t^{n-j+\alpha} e^{-t\rho(x,y)} \sum_{k=0}^{\infty} b_{j,k}(x,y) \, t^{-k},$$

with some $b_{j,k} \in C^\infty(\overline{\Omega} \times \overline{\Omega})$, where in particular

$$(41) \quad b_{j,0}(x,x) = t^{-\alpha} b_j(x,x)\sigma(P)(x,-t\nabla_x\rho(x,y)|_{y=x}) = b_j(x,x)\,\sigma(T_P)(x,\partial r(x)).$$

Thus (40) equals

$$(42) \qquad \int_0^\infty e^{-t\rho(\cdot,y)} \sum_{j=0}^{\infty} t^{n-j+\alpha} \, \widetilde{b}_j(\cdot,y) \, dt$$

with some $\widetilde{b}_j \in C^\infty(\overline{\Omega} \times \overline{\Omega})$, where

$$(43) \qquad \widetilde{b}_0(x,x) = b_0(x,x)\sigma(T_P)(x,\partial r(x)).$$

Combining this with the classical formulas (valid for $\mathrm{Re}\, p > 0$)

$$(44)$$
$$\mathrm{p.\,f.} \int_0^\infty e^{-tp} \, t^s \, dt = \begin{cases} \dfrac{\Gamma(s+1)}{p^{s+1}} & \text{if } s \in \mathbf{C}, \ s \neq -1,-2,\dots, \\[2ex] \dfrac{(-1)^{k+1}}{k!} \, p^k \, (\log p + C_k), & s = -k-1, \ k = 0,1,2,\dots \end{cases}$$

(where $C_k$ is a constant: $C_k = \lim_{m\to\infty} \sum_{j=k+1}^{m} \frac{1}{j} - \log m$), we thus obtain

$$(45) \quad K_w(x,y) = \begin{cases} \dfrac{a(x,y)}{\rho(x,y)^{n+\alpha+1}} + b(x,y) & \text{if } n+\alpha \notin \mathbf{Z}, \\[2ex] \dfrac{a(x,y)}{\rho(x,y)^{n+\alpha+1}} + b(x,y)\log\rho(x,y) & \text{if } n+\alpha \in \mathbf{Z}_{\geq 0}, \\[2ex] \dfrac{a(x,y)}{\rho(x,y)^{n+\alpha+1}} \log\rho(x,y) + b(x,y) & \text{if } n+\alpha \in \mathbf{Z}_{<0}, \end{cases}$$

which gives (10). (Actually, since in our case $n + \alpha + 1 > n > 0$, the third case in (45) never occurs.)

It remains to compute the boundary value of the leading term. This corresponds to the term $j = 0$ in (42), i.e., by (43) and (44),

$$\int_0^\infty e^{-t\rho(x,x)} \frac{t^{n+\alpha} b_0(x,x)}{\sigma(T_Q)(x, \partial r(x))} \, dt = \frac{\Gamma(n+\alpha+1)}{\rho(x,x)^{n+\alpha+1}} \frac{b_0(x,x)}{\sigma(T_Q)(x, \partial r(x))}$$

(with the first term on the right modified accordingly if $n + \alpha$ is a negative integer). Taking in particular $\alpha = 0$ and $g \equiv 0$, so that $Q = I$, and comparing with (2) shows that for $x \in \partial\Omega$,

$$b_0(x,x) = \frac{J[\rho](x)}{\pi^n}.$$

Consequently, for $x \in \partial\Omega$,

$$a(x,x) = \begin{cases} \dfrac{\Gamma(n+\alpha+1) \, J[\rho](x)}{\pi^n \, \sigma(T_Q)(x, \partial r(x))}, & n + \alpha \notin \mathbf{Z}_{<0}, \\[2ex] \dfrac{(-1)^{k+1}}{k!} \dfrac{J[\rho](x)}{\pi^n \sigma(T_Q)(x, \partial r(x))}, & n + \alpha = -k - 1, \, k = 0, 1, 2, \dots . \end{cases} \tag{46}$$

Substituting from (38) for $\sigma(T_Q)$, we get (11). (Again, the second case in (46) is not now needed.) This completes the proof. $\square$

*Remark 6.* As was already mentioned in the Introduction, the "off-diagonal part" of (10) — i.e. that $K_w$ is $C^\infty$ on $\overline{\Omega} \times \overline{\Omega}$ minus the boundary diagonal — was proved by Peloso [29], while the part concerning the singularity on the boundary diagonal was described by Komatsu (unpublished, referred to as "personal communication" in [16]). (It should be noted that in [16], the $C^\infty$ part in (10) for $\alpha \notin \mathbf{Z}$ was omitted by mistake.) $\square$

*Remark 7.* For $n = 1$, it is known that for the ordinary Bergman kernel (i.e. $\alpha = 0$, $g \equiv 0$) on a smoothly bounded domain $\Omega \subset \mathbf{C}$ the log term in (10) is actually absent. It should be noted, however, that in the case of *weighted* Bergman kernels the log-term appears even in dimension $n = 1$: for instance, for $\Omega = \mathbf{D}$, the unit disc, and $w(z) = 2 - |z|^2$, a short computation reveals that

$$K_w = \rho^{-2} - \rho^{-1} - 2\log\rho + C(\overline{\Omega}) \qquad (\rho(x,y) = 1 - x\overline{y}).$$

Similarly, for $m$ a positive integer it can be shown (using properties of hypergeometric functions) that $K_w$ for

$$w(z) = (1 - |z|^2)^m (2 - |z|^2)$$

on $\mathbf{D}$ contains a log term. $\square$

## 5. Sobolev-Bergman kernels

We now treat the reproducing kernels of the holomorphic Sobolev spaces $W^s_{\text{hol}}(\Omega)$ with respect to the norm (15) (for some nonnegative integer $m > 2s - 1$), which we denote by $\| \cdot \|^{\#}_{m,s}$.

To minimize confusion, for the rest of the paper we will reserve the symbol $\rho$ for the almost-sesquianalytic extension $\rho(x,y)$, denoting the single-variable function $\rho(x) = \rho(x,x) = -r(x)$ by $|r|$.

**Theorem 8.** *Let $W_{\mathrm{hol}}^s(\Omega)$ be the $s$-th holomorphic Sobolev space equipped with the norm*

$$\|f\|_{m,s}^{\#} := \Big[ \sum_{|\nu| \le m} \binom{|\nu|}{\nu} \|\partial^\nu f\|_{L^2(\Omega,|r|^{2m-2s})}^2 \Big]^{1/2},$$

*where $m$ is a nonnegative integer $> 2s - 1$ and $|\nu| = \nu_1 + \cdots + \nu_n$, $\binom{|\nu|}{\nu} = \frac{|\nu|!}{\nu_1! \ldots \nu_n!}$. Then the corresponding reproducing kernel $K^{(s)}(x,y)$ has the form (45) for $\alpha = -2s$, with some almost-sesquianalytic functions $a, b$ on $\overline{\Omega} \times \overline{\Omega}$ such that for $x \in \partial\Omega$*

$$(47) \qquad a(x,x) = \frac{\Gamma(n-2s+1)}{\Gamma(2m-2s+1)} \frac{J[\rho](x)}{\pi^n \|\partial r\|^{2m}},$$

*where $\Gamma(n-2s+1)$ is to be replaced by $\frac{(-1)^{k+1}}{k!}$ if $n-2s+1 = -k$, $k = 0, 1, 2, \ldots$.*

*Proof.* For any $f \in C_{\mathrm{hol}}^\infty(\overline{\Omega})$,

$$\|f\|_{m,s}^{\#2} = \sum_{|\nu| \le m} \binom{|\nu|}{\nu} \int_\Omega |\partial^\nu f|^2 \, |r|^{2m-2s}.$$

Integrating by parts, since $f$ is holomorphic, the integral equals to

$$(48) \qquad \int_{\partial\Omega} |r|^{2m-2s} \, \partial^\nu f \, \overline{\partial^{\nu-e_j} f} \, \overline{n}_j - \int_\Omega \overline{\partial}_j(|r|^{2m-2s}) \, \partial^\nu f \, \overline{\partial^{\nu-e_j} f},$$

where $j$ is any index for which $\nu_j \ge 1$, $e_j$ is the multiindex $(0, 0, \ldots, 1, \ldots, 0)$ with 1 in the $j$-th slot, and $n_j = \partial_j r / \|\partial r\|$ is the appropriate component of the outward unit normal. Since $2m - 2s > m - 1$, the first integral vanishes if $m \ge 1$, while the second can again be integrated by parts. Continuing in this fashion, we obtain after $|\nu|$ steps

$$(-1)^{|\nu|} \int_\Omega (\overline{\partial}^\nu |r|^{2m-2s}) \, \partial^\nu f \, \overline{f}.$$

We thus see that

$$(49) \qquad \|f\|_{m,s}^{\#2} = \langle \Theta f, f \rangle_{L^2(\Omega)} = \langle \boldsymbol{\Pi}\Theta f, f \rangle_{L_{\mathrm{hol}}^2(\Omega)}, \qquad \forall f \in C_{\mathrm{hol}}^\infty(\overline{\Omega}),$$

where

$$\Theta = \sum_{|\nu| \le m} \binom{|\nu|}{\nu} (-1)^{|\nu|} (\overline{\partial}^\nu |r|^{2m-2s}) \, \partial^\nu.$$

Using the restriction and extension operators $\gamma$ and $\mathbf{K}$ from the preceding section, the operator $\boldsymbol{\Pi}\Theta$ can again be viewed as an operator on $C_{\mathrm{hol}}^\infty(\partial\Omega)$. Namely, there exist tangential operators $Z_k$, $k = 1, \ldots, n$, such that

$$(50) \qquad \gamma\partial_k f = Z_k \gamma f, \qquad \forall f \in C_{\mathrm{hol}}^\infty(\overline{\Omega})$$

(or $\partial_k \mathbf{K} u = \mathbf{K} Z_k u \ \forall u \in C_{\mathrm{hol}}^\infty(\partial\Omega)$). Explicitly, one has

$$Z_k = \partial_k - \sum_{j=1}^n \frac{r_j r_k}{\|\partial r\|^2} \, \overline{\partial}_j,$$

where for brevity we have introduced the notation

$$r_j := \partial_j r.$$

The operator $\mathbf{\Pi}\Theta$, viewed as an operator on $C^\infty_{\mathrm{hol}}(\partial\Omega)$, is then

$$\gamma\mathbf{\Pi}\Theta\mathbf{K} = \sum_{|\nu|\leq m} \binom{|\nu|}{\nu} (-1)^{|\nu|} \gamma\mathbf{\Pi}(\overline{\partial}^\nu |r|^{2m-2s})\mathbf{K}Z^\nu$$

$$= \sum_{|\nu|\leq m} \binom{|\nu|}{\nu} (-1)^{|\nu|} \gamma\mathbf{T}_{\overline{\partial}^\nu |r|^{2m-2s}}\mathbf{K}Z^\nu.$$

(Here $Z^\nu$, of course, stands for $Z_1^{\nu_1}\ldots Z_n^{\nu_n}$, and similarly for $(\partial r)^\nu$ below.) Now

$$\overline{\partial}^\nu |r|^{2m-2s} = (-1)^{|\nu|}|r|^{2m-|\nu|-2s}\frac{\Gamma(2m-2s+1)}{\Gamma(2m-|\nu|-2s+1)}[(\overline{\partial r})^\nu + r\, g_\nu],$$

where $g_\nu \in C^\infty(\overline{\Omega})$. By the results reviewed in the preceding section, we thus see that $\gamma\mathbf{T}_{\overline{\partial}^\nu |r|^{2m-2s}}\mathbf{K}|_{H^2(\partial\Omega)}$ is, modulo smoothing operators, a generalized Toeplitz operator on $H^2(\partial\Omega)$ of order $2s-2m+|\nu|$ with principal symbol

$$\Gamma(2m-2s+1)\|\partial r\|^{2m-|\nu|-2s}\|\xi\|^{2s-2m+|\nu|}(-1)^{|\nu|}(\overline{\partial r})^\nu.$$

Since $Z^\nu$ is a differential operator of order $|\nu|$, it follows that

$$\gamma\mathbf{\Pi}\Theta\mathbf{K}|_{H^2(\partial\Omega)} \sim T_Q,$$

where $Q$ is a $\Psi$DO on $\partial\Omega$ of order $2s$ and with symbol satisfying

$$\sigma(Q)|_\Sigma = \Gamma(2m-2s+1) \sum_{|\nu|=m} \binom{m}{\nu} \frac{\|\partial r\|^{m-2s}}{\|\xi\|^{m-2s}}(\overline{\partial r})^\nu \cdot \sigma(Z^\nu)|_\Sigma.$$

However, for $(x,\xi) \in \Sigma$, i.e. $x \in \Omega$ and $\xi = t\eta_x$, $t > 0$, we have

$$\sigma(Z_k)(x,\xi) = \langle t\eta_x, Z_k\rangle = tr_k = \frac{\|\xi\|}{\|\partial r\|}r_k.$$

Thus

$$\sigma(T_Q) = \Gamma(2m-2s+1)\|\partial r\|^{m-2s}\|\xi\|^{2s}\sum_{|\nu|=m}\binom{m}{\nu}(\overline{\partial r})^\nu\frac{(\partial r)^\nu}{\|\partial r\|^m}$$

$$(51) \qquad = \Gamma(2m-2s+1)\|\xi\|^{2s}\|\partial r\|^{2m-2s}.$$

It follows that $T_Q$ is elliptic and has a parametrix $T_P$ of order $-2s$; by (P2) we may again assume that $P\Pi = \Pi P$, so that $T_P$ is just the restriction of $P$ to holomorphic functions.

From the ellipticity it follows that $T_Q$ and, hence, also $\mathbf{\Pi\Theta}|_{L^2_{\text{hol}}(\Omega)}$ is Fredholm as an operator from $W^{2s}_{\text{hol}}(\Omega)$ into $L^2_{\text{hol}}(\Omega)$. On the other hand, since by (49) $\mathbf{\Pi\Theta}|_{L^2_{\text{hol}}(\Omega)} \geq \mathbf{T}_{|r|^{-2s}} > 0$, the symmetric operator $\mathbf{\Pi\Theta}|_{L^2_{\text{hol}}(\Omega)}$ is injective and has self-adjoint closure, hence also dense range, as an operator on $L^2_{\text{hol}}(\Omega)$. Consequently, $\mathbf{\Pi\Theta}$ maps $W^{2s}_{\text{hol}}(\Omega)$ onto $L^2_{\text{hol}}(\Omega)$. Applying now part (iii) of Proposition 2 with $\Lambda = (\mathbf{\Pi\Theta\Pi})^{1/2}$, so that $\gamma\Lambda^2\mathbf{K}|_{H^2(\partial\Omega)} = \gamma\mathbf{\Pi\Theta K}|_{H^2(\partial\Omega)} \sim T_Q$, we get as before

$$K^{(s)}_y := K^{(s)}(\,\cdot\,,y) \sim PK_y,$$

and the same argument as in the proof of Theorem A (using the formulas (44) and (45)) yields the desired conclusion. The leading term (47) is also evaluated in the same manner, using (46) and (51).  $\square$

We next address the norm given by (16), which we denote by $\|\cdot\|^\flat_{m,s}$.

**Theorem 9.** *Let $W^s_{\text{hol}}(\Omega)$ be the $s$-th holomorphic Sobolev space equipped with the norm*

$$\|f\|^\flat_{m,s} := \Big[ \sum_{j=0}^m \|\mathcal{D}^j f\|^2_{L^2(\Omega,|r|^{2m-2s})} \Big]^{1/2},$$

*where $m$ is a nonnegative integer $> 2s-1$ and $\mathcal{D}$ is the "normal derivative" operator*

$$\mathcal{D} = \sum_{j=1}^n \overline{r}_j\,\partial_j.$$

*Then the corresponding reproducing kernel $K^{(s)}(x,y)$ has the form (45) for $\alpha = -2s$, with some almost-sesquianalytic functions $a,b$ on $\overline\Omega \times \overline\Omega$ such that for $x \in \partial\Omega$*

$$(52)\qquad a(x,x) = \frac{\Gamma(n-2s+1)}{\Gamma(2m-2s+1)}\,\frac{J[\rho](x)}{\pi^n\|\partial r\|^{4m}},$$

*where $\Gamma(n-2s+1)$ is to be replaced by $\frac{(-1)^{k+1}}{k!}$ if $n-2s+1 = -k$, $k = 0,1,2,\dots$.*

*Proof.* The argument is similar to the one for the preceding theorem, so we will be brief. By partial integration, we have for any $f \in C^\infty_{\text{hol}}(\overline\Omega)$,

$$\int_\Omega |\mathcal{D}^j f|^2\,|r|^{2m-2s} = \int_{\partial\Omega} |r|^{2m-2s}\,\mathcal{D}^j f\,\overline{\mathcal{D}^{j-1}f}\,\overline{n}_\mathcal{D} - \int_\Omega \mathcal{D}^*(|r|^{2m-2s}\mathcal{D}^j f)\,\overline{\mathcal{D}^{j-1}f},$$

where $n_\mathcal{D} = \frac{\mathcal{D}r}{\|\partial r\|} = \|\partial r\|$ is the appropriate component of the outward unit normal, and $\mathcal{D}^*$ is the formal adjoint

$$\mathcal{D}^*g := \sum_{k=1}^n \overline{\partial}_k(r_k g) = \frac{\Delta r}{4}\,g + \sum_{k=1}^n r_k\,\overline{\partial}_k g = \Big(\frac{\Delta r}{4} + \overline{\mathcal{D}}\Big)g.$$

Since $m > 2s-1$ by hypothesis, the integral over $\partial\Omega$ vanishes if $m \geq 1$, while the second integral can again be integrated by parts. Continuing in this fashion, we obtain after $j$ steps

$$(-1)^j \int_\Omega \mathcal{D}^{*j}(|r|^{2m-2s}\mathcal{D}^j f)\,\overline{f}.$$

We thus again see that

$$\|f\|_{m,s}^{b2} = \langle \mathbf{\Pi}\Theta f, f\rangle_{L^2_{\mathrm{hol}}(\Omega)} \qquad \forall f \in C^\infty_{\mathrm{hol}}(\overline{\Omega}),$$

where $\Theta$ is given by

$$\Theta f = \sum_{j=0}^m (-1)^j \mathcal{D}^{*j}(|r|^{2m-2s}\mathcal{D}^j f).$$

Observe that for any $g \in C^\infty(\overline{\Omega})$ and any real number $a$,

$$\mathcal{D}^*(|r|^a g) = |r|^a \mathcal{D}^* g + g\overline{\mathcal{D}}|r|^a = |r|^a \mathcal{D}^* g - a|r|^{a-1}\|\partial r\|^2 g,$$

since $\overline{\mathcal{D}}|r|^a = -a|r|^{a-1}\overline{\mathcal{D}}r = -a|r|^{a-1}\|\partial r\|^2$. By a straightforward induction argument,

$$\mathcal{D}^{*j}(|r|^{2m-2s}\mathcal{D}^j f) = (-1)^j \frac{\Gamma(2m-2s+1)}{\Gamma(2m-j-2s+1)}\,|r|^{2m-j-2s}\,\|\partial r\|^{2j}\,\mathcal{D}^j f$$
$$+ |r|^{2m-j-2s+1}\mathcal{L}_{m,j}\mathcal{D}^j f,$$

where $\mathcal{L}_{m,j}$ is a differential operator with $C^\infty(\overline{\Omega})$ coefficients and involving only anti-holomorphic derivatives. In a similar fashion it can be shown that

$$\mathcal{D}^j f = \sum_{|\nu|=j}\binom{j}{\nu}(\overline{\partial r})^\nu\,\partial^\nu f + \sum_{|\nu|<j} a_{j,\nu}\partial^\nu f,$$

with some coefficients $a_{j,\nu} \in C^\infty(\overline{\Omega})$ independent of $f$. Combining everything together and using (50) we thus get

$$\gamma\mathbf{\Pi}\Theta\mathbf{K}|_{L^2_{\mathrm{hol}}(\Omega)} = \sum_{j=0}^m \frac{\Gamma(2m-2s+1)}{\Gamma(2m-j-2s+1)}\gamma\mathbf{\Pi}|r|^{2m-j-2s}\|\partial r\|^{2j}\sum_{|\nu|=j}\binom{j}{\nu}(\overline{\partial r})^\nu\mathbf{K}Z^\nu$$
$$+ \sum_{j=0}^m \frac{\Gamma(2m-2s+1)}{\Gamma(2m-j-2s+1)}\gamma\mathbf{\Pi}|r|^{2m-j-2s}\|\partial r\|^{2j}\sum_{|\nu|<j} a_{j,\nu}\mathbf{K}Z^\nu$$
$$+ \sum_{j=0}^m (-1)^j\gamma\mathbf{\Pi}|r|^{2m-j-2s+1}\sum_{|\nu|\le j} b_{j,\nu}\mathbf{K}Z^\nu,$$

with $b_{j,\nu} \in C^\infty(\overline{\Omega})$. As in the preceding proof, the second and the third term, as well as the summands $j = 0, \ldots, m-1$ of the first term, are (modulo smoothing errors) generalized Toeplitz operators of orders $\le 2s-1$, while the $j = m$ summand of the first term is a generalized Toeplitz operator of order $2s$ with symbol

$$\Gamma(2m-2s+1)\frac{\|\partial r\|^{m-2s}}{\|\xi\|^{m-2s}}\|\partial r\|^{2m}\sum_{|\nu|=m}\binom{m}{\nu}(\overline{\partial r})^\nu\left(\frac{\|\xi\|}{\|\partial r\|}\partial r\right)^\nu$$

$$\tag{53} = \Gamma(2m-2s+1)\,\|\xi\|^{2s}\|\partial r\|^{4m-2s}.$$

Thus $\gamma\mathbf{\Pi}\Theta\mathbf{K}|_{H^2(\partial\Omega)} \sim T_Q$, where $Q$ is a $\Psi$DO on $\partial\Omega$ of order $2s$ with symbol (53). The same argument involving the parametrix and the computation of the leading term (using the formulas (45) and (46), respectively) complete the proof. $\quad\square$

## 6. SOME MORE KERNELS

We digress to make some observations which will not be needed in the sequel, but which yield an independent proof of the equivalence of the Sobolev norm on holomorphic functions with the norm (15) mentioned in the Introduction (though not of the equivalence with the norm (16) — we will, however, give an independent proof of that one too in the next section), as well as an extension of this equivalence to harmonic functions which seems to be new.

Recall that if $\Theta$ is a positive definite (or just "positive" for short — meaning that $\langle \Theta f, f \rangle > 0$ for any $f \neq 0$ in dom $\Theta$) self-adjoint operator on $L^2(\partial\Omega)$ which is an elliptic $\Psi$DO of order $m \neq 0$, then the square root $\Theta^{1/2}$ of $\Theta$ (in the sense of the spectral theorem) is also a $\Psi$DO, of order $m/2$. Since the positivity condition and self-adjointness imply that $\Theta$ has an inverse with the same properties, it also follows that the negative square root $\Theta^{-1/2}$ is a $\Psi$DO of order $-m/2$. These facts are actually just special cases of the more general result, going back to Seeley, asserting that the complex power $\Theta^s$ is a $\Psi$DO of order $m \operatorname{Re} s$, for any $s \in \mathbf{C}$; we will have more to say about this in the next section. Note that the ellipticity again implies that $\Theta$ is Fredholm as an operator from each $W^s(\partial\Omega)$, $s \in \mathbf{R}$, into $W^{s-m}(\partial\Omega)$; being injective and, hence, with dense range as an operator on $L^2(\partial\Omega)$, $\Theta$ is therefore a bijection of $W^m(\partial\Omega)$ onto $L^2(\partial\Omega)$. In particular, dom $\Theta = W^m(\partial\Omega)$.

**Proposition 10.** *Let $\Theta$ be a positive self-adjoint operator on $L^2(\partial\Omega)$ which is an elliptic $\Psi$DO of order $m \in \mathbf{R}$. Let $\mathcal{H}_\Theta$ be the completion of $C^\infty(\partial\Omega) \subset \operatorname{dom} \Theta$ with respect to the norm*

$$\|u\|_\Theta^2 := \langle \Theta u, u \rangle.$$

*Then $\mathcal{H}_\Theta = W^{m/2}(\partial\Omega)$, with equivalent norms.*

*Proof.* For $m = 0$, both $\Theta$ and, by ellipticity, $\Theta^{-1}$ are $\Psi$DOs of order 0, hence bounded, and it follows immediately that $\|\cdot\|_\Theta$ is equivalent to the norm in $L^2(\partial\Omega)$. We may thus assume that $m \neq 0$. Hence, by the remarks above, $\Theta^{1/2}$ is an elliptic $\Psi$DO of order $m/2$. For any $u \in C^\infty(\partial\Omega) \subset \operatorname{dom} \Theta \subset \operatorname{dom} \Theta^{1/2}$, we have $\|u\|_\Theta = \|\Theta^{1/2} u\|_{L^2(\partial\Omega)}$. By ellipticity, the last norm is equivalent to $\|u\|_{W^{m/2}(\partial\Omega)}$. The claim follows. $\square$

**Corollary 11.** *For $\Theta$ as in the preceding proposition, let $\mathcal{L}_\Theta$ be the completion of $C^\infty_{\mathrm{harm}}(\overline{\Omega})$, the subspace of harmonic functions in $C^\infty(\overline{\Omega})$, with respect to the norm*

$$\|f\|_\Theta^2 := \langle \Theta \gamma f, \gamma f \rangle_{L^2(\partial\Omega)}.$$

*Then $\mathcal{L}_\Theta = W^{(m+1)/2}_{\mathrm{harm}}(\Omega)$, with equivalent norms.*

*Proof.* The mapping $u \mapsto \mathbf{K}u$ (with the inverse $f \mapsto \gamma f$) is known to be an isomorphism of $W^s(\partial\Omega)$ onto $W^{s+\frac{1}{2}}(\Omega)$, for any $s \in \mathbf{R}$; see Lions and Magenes [27], Chapter 2, §7.3. $\square$

As an application, we immediately get the promised independent proof of the equivalence of the norms (15) with the Sobolev norms on holomorphic functions, as well as an analogous result for the harmonic functions.

In addition to the "normal derivative" $\mathcal{D} = \sum_j \overline{r}_j \partial_j$, let us introduce also the operator

$$\overline{\mathcal{D}} := \sum_{j=1}^n r_j \overline{\partial}_j.$$

**Theorem 12.** *Let $s$ be a real number, $m > s - \frac{1}{2}$ a nonnegative integer, and $x_0$ any point in $\Omega$. Then a harmonic function $f$ belongs to $W^s(\Omega)$ if and only if any of the following quantities is finite, and the square root of each of these quantities gives an equivalent norm in $W_{\mathrm{harm}}^s(\Omega)$:*

  (a) $\sum_{|\nu|+|\mu| \leq m} \frac{|\nu+\mu|!}{\nu!\mu!} \|\partial^\nu \overline{\partial}^\mu f\|_{L^2(\Omega,|r|^{2m-2s})}^2$;

  (b) $\sum_{|\nu|+|\mu|=m} \frac{m!}{\nu!\mu!} \|\partial^\nu \overline{\partial}^\mu f\|_{L^2(\Omega,|r|^{2m-2s})}^2 + \|f\|_{L^2(\Omega,|r|^{2m-2s})}^2$;

  (c) $\sum_{|\nu|+|\mu|=m} \frac{m!}{\nu!\mu!} \|\partial^\nu \overline{\partial}^\mu f\|_{L^2(\Omega,|r|^{2m-2s})}^2 + \sum_{|\nu|+|\mu|<m} |\partial^\nu \overline{\partial}^\mu f(x_0)|^2$;

  (d) $\sum_{j=0}^m \|(\mathcal{D} + \overline{\mathcal{D}})^j f\|_{L^2(\Omega,|r|^{2m-2s})}^2$;

  (e) $\|(\mathcal{D} + \overline{\mathcal{D}})^m f\|_{L^2(\Omega,|r|^{2m-2s})}^2 + \|f\|_{L^2(\Omega,|r|^{2m-2s})}^2$.

*Proof.* (a) Let us define operators $R_j$ and $\overline{R}_j$, $j = 1, \ldots, n$, on $L^2(\partial\Omega)$ by

$$R_j u = \gamma \partial_j \mathbf{K} u, \qquad \overline{R}_j u = \gamma \overline{\partial}_j \mathbf{K} u$$

(where $\gamma$, $\mathbf{K}$ and, below, $\Lambda_0 = \mathbf{K}^*\mathbf{K}$ have the same meaning as in Section 4). Setting as usual $f = \mathbf{K}u$, we have

$$(54) \quad \begin{aligned} \|\partial^\nu \overline{\partial}^\mu f\|_{L^2(\Omega,|r|^{2m-2s})}^2 &= \|\mathbf{K}R^\nu \overline{R}^\mu u\|_{L^2(\Omega,|r|^{2m-2s})}^2 = \||r|^{m-s}\mathbf{K}R^\nu \overline{R}^\mu u\|_{L^2(\Omega)}^2 \\ &= \langle \overline{R}^{*\mu} R^{*\nu} \mathbf{K}^* |r|^{2m-2s}\mathbf{K}R^\nu \overline{R}^\mu u, u\rangle_{L^2(\partial\Omega)}. \end{aligned}$$

Let us introduce the notation, for any expression $X$,

$$\mathcal{R}[X] := \sum_{j=1}^n (R_j^* X R_j + \overline{R}_j^* X \overline{R}_j).$$

By (54), the norm in (a) can then be written as $\langle \Theta u, u\rangle$, where

$$(55) \qquad \Theta = \sum_{j=0}^m \mathcal{R}^j[\mathbf{K}^* |r|^{2m-2s}\mathbf{K}].$$

As we have seen in Section 4, $\mathbf{K}^* |r|^{2m-2s}\mathbf{K}$ is a $\Psi$DO on $\partial\Omega$ of order $2s - 2m - 1$ with symbol $\Gamma(2m - 2s + 1)\|\xi\|^{2s-2m-1}\|\partial r\|^{2m-2s}$. On the other hand, $R_j$ and $\overline{R}_j$ are $\Psi$DOs of order 1, and by Green's theorem, for any $f \in C_{\mathrm{harm}}^\infty(\overline{\Omega})$

$$\int_\Omega \sum_{j=1}^n |\partial_j f|^2 = \int_{\partial\Omega} \sum_{j=1}^n \overline{n}_j\, \partial_j f\, \overline{f}$$

(where $\overline{n}_j = \overline{r}_j / \|\partial r\|$), or

$$\sum_{j=1}^{n} R_j^* \mathbf{K}^* \mathbf{K} R_j = \sum_{j=1}^{n} \overline{n}_j R_j.$$

Taking symbols gives

$$\sum_{j=1}^{n} |\sigma(R_j)|^2 \sigma(\Lambda_0) = \sum_{j=1}^{n} \overline{n}_j \sigma(R_j).$$

Similarly for $\overline{R}_j$. (Note that $\sigma(\overline{R}_j)$ is *not* the complex conjugate of $\sigma(R_j)$.) Consequently, for any $X$

$$\begin{aligned}
\sigma(\mathcal{R}[X]) &= \sum_{j=1}^{n} (|\sigma(R_j)|^2 + |\sigma(\overline{R}_j)|^2) \sigma(X) \\
&= \sum_{j=1}^{n} \frac{\overline{n}_j \sigma(R_j) + n_j \sigma(\overline{R}_j)}{2\sigma(\Lambda_0)} \, \sigma(X) \\
&= \frac{\sigma(\vartheta)}{4\sigma(\Lambda_0)} \, \sigma(X),
\end{aligned}$$

where $\vartheta = 2\sum_{j=1}^{n}(\overline{n}_j R_j + n_j \overline{R}_j) = \gamma \frac{\partial}{\partial n} \mathbf{K}$ is the classical Dirichlet-to-Neumann operator, which is an elliptic $\Psi$DO of order 1. Thus $\Theta$ is a $\Psi$DO of order $2s - 1$ with symbol

$$\begin{aligned}
\|\xi\|^m \sigma(\vartheta/2)^m \Gamma(2m - 2s + 1)\|\xi\|^{2s-2m-1}\|\partial r\|^{2m-2s} \\
= \Gamma(2m - 2s + 1)\|\xi\|^{2s-1}\Big(\frac{\sigma(\vartheta)}{\|2\xi\|}\Big)^m \|\partial r\|^{2m-2s} > 0,
\end{aligned}$$

so $\Theta$ is elliptic. At the same time, as an operator on $L^2(\partial\Omega)$, $\Theta$ is nonnegative self-adjoint (being a sum of products of the form $V^*V$, where $V : L^2(\partial\Omega) \to L^2(\Omega)$ given by $u \mapsto |r|^{m-s}\partial^\nu\overline{\partial}^\mu \mathbf{K}$ is densely-defined and closed) and positive (since the summand $\nu = \mu = 0$ gives just the norm of $f$ in $L^2(\Omega, |r|^{2m-2s})$). By Corollary 11, the claim follows.

(b) This is the same as in (a), except that we are keeping from $\Theta$ only the terms of the highest ($|\nu| + |\mu| = m$) and lowest ($|\nu| = |\mu| = 0$) order; since the latter was responsible for $\Theta$ being (not only nonnegative but) positive, while the former was responsible for $\Theta$ having the right order $2s - 1$ and elliptic symbol, the conclusion remains in force.

(c) This time $\Theta$ is of the form $\Theta = \Theta' + \Theta''$, where $\Theta'$ is again the top degree part of (55), while

$$\Theta'' = \sum_{|\nu|+|\mu|<m} \langle \, \cdot \, , P_{x_0}^{\nu\overline{\mu}} \rangle P_{x_0}^{\nu\overline{\mu}},$$

where $P_{x_0}^{\nu\overline{\mu}}(\zeta) := \partial^\nu\overline{\partial}^\mu P(\zeta, \cdot)|_{x_0}$ is the derivative at $x_0$ of the Poisson kernel $P(\zeta, x)$. Since $P$ is known to be $C^\infty$ on $\partial\Omega \times \Omega$, $\Theta''$ is a smoothing operator. Thus $\Theta$ is again an elliptic $\Psi$DO of order $2s - 1$, and since $\Theta'$ and $\Theta''$ are both self-adjoint and nonnegative, while $\mathrm{Ker}\,\Theta' = \{$polynomials of degree $< m\}$ whereas $\mathrm{Ker}\,\Theta'' = \{$functions vanishing at $x_0$ to order at least $m\}$, $\Theta = \Theta' + \Theta''$ is positive. Thus the claim again follows by Corollary 11.

(d) Here matters are more complicated since $(\mathcal{D} + \overline{\mathcal{D}})f$ is no longer harmonic when $f$ is. However, by the Leibniz rule, we have

$$(\mathcal{D} + \overline{\mathcal{D}})^j = \sum_{|\nu| + |\mu| \leq j} a_{j\nu\mu} \partial^\nu \overline{\partial}^\mu$$

where for $|\nu| + |\mu| = j$, $a_{j\nu\mu} = \frac{j!}{\nu!\mu!}(\overline{\partial r})^\nu (\partial r)^\mu$. Thus again, for $f \in C_{\mathrm{harm}}^\infty(\overline{\Omega})$,

$$\|(\mathcal{D} + \overline{\mathcal{D}})^j f\|_{L^2(\Omega, |r|^{2m-2s})}^2 = \left\| |r|^{m-s} \sum_{|\nu| + |\mu| \leq j} a_{j\nu\mu} \partial^\nu \overline{\partial}^\mu f \right\|_{L^2(\Omega)}^2$$

$$= \left\| |r|^{m-s} \sum_{|\nu| + |\mu| \leq j} a_{j\nu\mu} \mathbf{K} R^\nu \overline{R}^\mu u \right\|_{L^2(\partial\Omega)}^2 = \langle \Theta_j u, u \rangle_{L^2(\partial\Omega)},$$

where

$$\Theta_j = \sum_{|\nu| + |\mu| \leq j} \sum_{|\epsilon| + |\eta| \leq j} \overline{R}^{*\eta} R^{*\epsilon} \mathbf{K}^* \overline{a}_{j\epsilon\eta} |r|^{2m-2s} a_{j\nu\mu} \mathbf{K} R^\nu \overline{R}^\mu.$$

Each summand is a $\Psi$DO on $\partial\Omega$ of order $\leq |\mu| + |\nu| + |\epsilon| + |\eta| + 2s - 2m - 1 \leq 2j + 2s - 2m - 1 \leq 2s - 1$, with equality occurring only for $|\epsilon| + |\eta| = |\mu| + |\nu| = j = m$; the corresponding symbol in that case is

$$\sigma(\Theta_m) = \sum_{|\nu| + |\mu| = m} \sum_{|\epsilon| + |\eta| = m} \overline{\sigma(\overline{R})^\eta \sigma(R)^\epsilon} \sigma(\mathbf{K}^* \overline{a}_{m\epsilon\eta}) |r|^{2m-2s} a_{m\nu\mu} \mathbf{K}) \sigma(R)^\nu \sigma(\overline{R})^\mu$$

$$= \Gamma(2m - 2s + 1) \|\xi\|^{2s-2m-1} \|\partial r\|^{2m-2s} \left| \sum_{|\nu| + |\mu| = m} a_{m\nu\mu} \sigma(R)^\nu \sigma(\overline{R})^\mu \right|^2$$

$$= \Gamma(2m - 2s + 1) \|\xi\|^{2s-2m-1} \|\partial r\|^{2m-2s}$$

$$\cdot \left| \sum_{|\nu| + |\mu| = m} \frac{m!}{\nu!\mu!} \sigma(\overline{n} \|\partial r\| R)^\nu \sigma(n \|\partial r\| \overline{R})^\mu \right|^2$$

$$= \Gamma(2m - 2s + 1) \|\xi\|^{2s-2m-1} \|\partial r\|^{2m-2s} \left| \|\partial r\| \sum_{j=1}^n \sigma(\overline{n}_j R_j + n_j \overline{R}_j) \right|^{2m}$$

$$= \Gamma(2m - 2s + 1) \|\partial r\|^{4m-2s} \left( \frac{\sigma(\vartheta)}{\|2\xi\|} \right)^{2m} \|\xi\|^{2s-1} > 0.$$

We thus see that the expression in (d) is of the form $\langle \Theta u, u \rangle_{L^2(\partial\Omega)}$ with $\Theta$ an elliptic $\Psi$DO of order $2s - 1$; further, $\Theta$ is again self-adjoint and nonnegative (for the same reason as in the proof of part (a)), and since the term $j = 0$ is just the norm of $f$ in $L^2(\Omega, |r|^{2m-2s})$, it is even positive. An application of Corollary 11 hence again yields the desired conclusion.

(e) This again follows from the proof of (d) for the same reason as (b) followed from the proof of (a): it is enough to keep the top degree term and the zero degree term from $\Theta$, since the former takes care of the right order and ellipticity, while the latter of injectivity and, hence, positivity. $\square$

*Remark 13.* Note that the operator $\mathcal{D} + \overline{\mathcal{D}}$ in (d) and (e) is essentially the usual (real) "normal derivative": namely, $(\mathcal{D} + \overline{\mathcal{D}})f = \frac{1}{2}\langle \nabla f, \nabla r \rangle_{\mathbf{R}^{2n}}$. $\square$

**Corollary 14.** *For $s$ a real number, $m > s - \frac{1}{2}$ a nonnegative integer, and $x_0$ any point of $\Omega$, a holomorphic function $f$ on $\Omega$ belongs to $W_{\mathrm{hol}}^s(\Omega)$ if and only if any of the quantities below is finite, and the square roots of these quantities are equivalent norms on $W_{\mathrm{hol}}^s(\Omega)$:*

(a) $\sum_{|\nu| \le m} \frac{|\nu|!}{\nu!} \|\partial^\nu f\|^2_{L^2(\Omega, |r|^{2m-2s})}$;

(b) $\sum_{|\nu| = m} \frac{m!}{\nu!} \|\partial^\nu f\|^2_{L^2(\Omega, |r|^{2m-2s})} + \|f\|^2_{L^2(\Omega, |r|^{2m-2s})}$;

(c) $\sum_{|\nu| = m} \frac{m!}{\nu!} \|\partial^\nu f\|^2_{L^2(\Omega, |r|^{2m-2s})} + \sum_{|\nu| < m} |\partial^\nu f(x_0)|^2$;

*and the quantities* (d) *and* (e) *from Theorem 12.*

*Proof.* Just specialize the preceding theorem to holomorphic functions. $\square$

Unfortunately, we have not been able to prove in this way the equivalence of (16) with the norm in $W_{\mathrm{hol}}^s(\Omega)$), though we expect this should be possible. (We will give a proof in the next section nonetheless, using generalized Toeplitz operators.)

With Theorem 12 in hands, the reproducing kernels can again be handled in the usual way.

**Theorem 15.** *The reproducing kernel of $W_{\mathrm{hol}}^s(\Omega)$ with respect to any of the norms in Corollary 14 is of the form* (45), *with the leading term given by* (46), *for $\alpha = -2s$ and $Q = \Lambda_0^{-1}\Theta$ where $\Theta$ is the operator from the proofs of* (a)–(e) *in Theorem 12.*

*Proof.* We have seen that the norms in question are of the form $\langle \Theta u, u \rangle_{L^2(\partial\Omega)}$, where $u = \gamma f$. In terms of $f = \mathbf{K}u$, this equals $\langle \Theta\gamma f, \gamma f \rangle_{L^2(\partial\Omega)} = \langle \gamma^*\Theta\gamma f, f \rangle_{L^2(\Omega)} = \langle \mathbf{K}\Lambda_0^{-1}\Theta\gamma f, f \rangle_{L^2(\Omega)} = \langle \mathbf{\Pi}\mathbf{K}\Lambda_0^{-1}\Theta\gamma f, f \rangle_{L^2_{\mathrm{hol}}(\Omega)}$. As before, viewed as an operator on $H^2(\partial\Omega)$, the last operator becomes

$$\gamma\mathbf{\Pi}\mathbf{K}\Lambda_0^{-1}\Theta\gamma\mathbf{K}|_{H^2(\partial\Omega)} = \Pi_{\Lambda_0}\Lambda_0^{-1}\Theta|_{H^2(\partial\Omega)},$$

which we have seen in Section 4 to be $\sim T_Q$ for $Q$ a $\Psi$DO having the same order and symbol on $\Sigma$ as $\Lambda_0^{-1}\Theta$. The assertion now follows in the same way as in the proof of Theorem A, using the formulas (45) and (46). $\square$

## 7. HOLOMORPHIC CONTINUATION

We finally deal with the kernels for the Sobolev-Bergman spaces with norms modelled on (13). This construction further makes sense even for complex $s$, thus yielding an analytic continuation of the kernels with respect to the parameter $s$.

Let $A$ be a positive self-adjoint elliptic $\Psi$DO of degree $m > 0$ on $\partial\Omega$. Then $A^{-1}$ is compact, hence its spectrum consists of isolated eigenvalues $0 < \lambda_1 < \lambda_2 < \ldots$ of finite multiplicity. We can therefore define for any $s \in \mathbf{C}$ the operator $A^s$ by the spectral theorem, i.e.

$$A^s = \sum_j \lambda_j^s P_j$$

where $P_j$ is the projection onto the eigenspace corresponding to $\lambda_j$. Alternatively, one can define $A^s$ for $\operatorname{Re} s < 0$ by the contour integral

$$A^s = \oint_{\lambda_1/2-i\infty}^{\lambda_1/2+i\infty} \lambda^s (A-\lambda)^{-1} \, d\lambda$$

(with the branch of $\lambda^s$ defined in the right half-plane so that $1^s = 1$). For $\operatorname{Re} s \geq 0$, one then sets

$$A^s = A^k A^{s-k}, \qquad k > \operatorname{Re} s;$$

this is unambiguous since $A^s A = A^{s+1}$ for $\operatorname{Re} s < -1$.

For a positive self-adjoint elliptic $\Psi$DO of degree $m < 0$, one then defines $A^s$ as $(A^{-1})^{-s}$, the right-hand side being taken in the sense of the previous paragraph. In both cases ($m < 0$ and $m > 0$), the operator $A^s$ so defined is normal for any $s \in \mathbf{C}$, and self-adjoint and positive if $s$ is real.

It is then a result going back to Seeley [33], with some later developments, stated in a form convenient for our purpose, in Bucicovschi [9] or Schrohe [32], that the operator $A^s$ defined as above is again a $\Psi$DO, of order $ms$, and with symbol $\sigma(A)^s$.

It turns out that all this remains true also for generalized Toeplitz operators.

**Proposition 16.** *Let $T$ be a positive self-adjoint operator on $H^2(\partial\Omega)$ such that $T \sim T_Q$, where $T_Q$ is of degree $m \neq 0$ and elliptic with $\sigma(T_Q) > 0$. Let $T^s$, $s \in \mathbf{C}$, be defined using the spectral theorem. Then $T^s \sim T_{Q_s}$, where $Q_s$ is an elliptic $\Psi$DO on $\partial\Omega$ of order $ms$, and $\sigma(T_{Q_s}) = \sigma(T_Q)^s$.*

*Proof.* Replacing $Q$ by $(Q + Q^*)/2$, we can assume that $Q$ is self-adjoint. Since $T_Q \sim T_{Q'}$ if the total symbols of $Q$ and $Q'$ agree in a neighbourhood of $\Sigma$, we may also assume that $\sigma(Q) > 0$ not only on $\Sigma$ but everywhere, i.e. that $Q$ is elliptic. Finally, arguing as on p. 20 in [7] we can assume that $Q$ also commutes with $\Pi$. (More precisely: it is obvious from the microlocal model that one can always choose microlocally $Q_j$ with positive symbols and a $\Psi$DO partition of unity $\chi_j$ which commute with $\Pi$ modulo smoothing operators. Then $Q = \sum \chi_j Q_j$ has a positive symbol, and $[Q, \Pi] \sim 0$; replacing $Q$ by $\Pi Q \Pi + (I-\Pi)Q(I-\Pi) = Q + (2\Pi-1)[Q,\Pi]$, which is a $\Psi$DO $\sim Q$, we can thus assume that $Q\Pi = \Pi Q$ and $\sigma(Q) > 0$. The author is grateful to Louis Boutet de Monvel for this explanation.) By the spectral theorem, the powers $(Q^2)^{s/2}$ (of the positive selfadjoint operator $Q^2$) then also commute with $\Pi$, and $\Pi(Q^2)^{s/2}\Pi = (\Pi Q^2\Pi)^{s/2} = (\Pi Q\Pi)^s$, for any $s \in \mathbf{C}$. Setting $Q_s := (Q^2)^{s/2}$, the result follows. $\square$

**Theorem 17.** *Let $T$ be a positive self-adjoint operator on $H^2(\partial\Omega)$ such that $T \sim T_P$, where $\sigma(T_P) > 0$ and $\operatorname{ord}(T_P) = 2s - 1$, $s \in \mathbf{R}$. Let $\mathcal{H}_T$ be the completion of $C^\infty_{\mathrm{hol}}(\overline{\Omega})$ with respect to the norm*

$$\|f\|_T^2 := \langle T\gamma f, \gamma f\rangle_{H^2(\partial\Omega)}. \tag{56}$$

*Then*

(a) *$\mathcal{H}_T$ coincides with the holomorphic Sobolev space $W^s_{\mathrm{hol}}(\Omega)$, with equivalent norms;*

(b) *the reproducing kernel of $\mathcal{H}_T$ with respect to the norm (56) has the form (45), with the leading term given by (46), for $\alpha = -2s$ and $Q = \Lambda_0^{-1}P$.*

*Proof.* In terms of $u = \gamma f \in C^\infty_{hol}(\partial\Omega)$, (56) becomes

$$(57) \qquad \|\mathbf{K}u\|^2_T = \langle Tu, u \rangle_{H^2(\partial\Omega)} = \|T^{1/2}u\|^2_{H^2(\partial\Omega)}.$$

By the preceding proposition, $T^{1/2} \sim T_{P_{1/2}}$ where $P_{1/2}$ is of order $s - \frac{1}{2}$ and $\sigma(T_{P_{1/2}}) > 0$. By (P7), the generalized Toeplitz operator $T_{P_{1/2}}$ is elliptic, and, hence, Fredholm as an operator from $W^{s-1/2}_{hol}(\partial\Omega)$ into $H^2(\partial\Omega)$, i.e. its range is closed and its kernel and cokernel are finite-dimensional. The same is therefore true for its compact perturbation $T^{1/2}$; since we know $T$ and, hence, $T^{1/2}$ to be positive (hence, injective) and self-adjoint (hence, by injectivity, with dense range) as an operator on $H^2(\partial\Omega)$, it follows that $T^{1/2}$ is an isomorphism of $W^{s-1/2}_{hol}(\Omega)$ onto $H^2(\partial\Omega)$. By (57), the space $\mathcal{H}_T$ (when its elements are viewed as functions on $\partial\Omega$ rather than $\Omega$, via the boundary values) thus coincides with $W^{s-\frac{1}{2}}_{hol}(\partial\Omega)$, with equivalent norms. Since the Poisson extension operator $\mathbf{K}$ is an isomorphism of $W^{s-\frac{1}{2}}_{harm}(\partial\Omega)$ onto $W^s_{hol}(\Omega)$, for any $s \in \mathbf{R}$ (this is immediate from the corresponding fact for $W^{s-1/2}_{harm}(\partial\Omega)$ and $W^s_{harm}(\Omega)$, recalled already in the proof of Corollary 11), the first part of the theorem follows.

For the second part, write

$$\|f\|^2_T = \langle \gamma^* T \gamma f, f \rangle_{L^2(\Omega)} = \langle \mathbf{K}\Lambda_0^{-1} T \gamma f, f \rangle_{L^2(\Omega)}$$
$$= \langle \mathbf{\Pi K}\Lambda_0^{-1} T \gamma f, f \rangle_{L^2_{hol}(\Omega)},$$

and note that the operator $\mathbf{\Pi K}\Lambda_0^{-1} T \gamma$, viewed as an operator on $H^2(\partial\Omega)$, becomes

$$\gamma \mathbf{\Pi K}\Lambda_0^{-1} T \gamma \mathbf{K} = \gamma \mathbf{\Pi K}\Lambda_0^{-1} T = \Pi_{\Lambda_0}\Lambda_0^{-1} T$$
$$\sim \Pi_{\Lambda_0}\Lambda_0^{-1} T_Q \sim T^{(\Lambda_0)}_{\Lambda_0^{-1}Q} \sim T_{\Lambda_0^{-1}Q}$$

(where, in the penultimate equivalence, we again used (P2) and assumed without loss of generality that $\Pi Q = Q\Pi$). Part (b) of the theorem therefore follows in the usual way from the formulas (45) and (46). $\quad\square$

*Remark 18.* Alternatively, the second part of the theorem could be proved directly — without passing from $H^2(\partial\Omega)$ to $L^2_{hol}(\Omega)$ — by working in $H^2(\partial\Omega)$ and using the Szegö kernel instead of the Bergman kernel: that is, using the analogue

$$K^{(T)}_y = T^{-1}S_y$$

of the formula (19) for the Szegö kernel $S_y(x) \equiv S(x, y)$ (whose proof is the same as the proof of (19)) and the corresponding analogue

$$S(x, y) = \int_0^\infty e^{-t\rho(x,y)} \sum_{j=0}^\infty t^{n-1-j} \, b^{(S)}_j(x, y)\, dt, \qquad b^{(S)}_j \in C^\infty(\overline{\Omega} \times \overline{\Omega}),$$

of the formula (39), which was also proved in [8]. $\quad\square$

As a first application of Theorem 17, we have the following corollary. (Parts (a)–(c) are the same as in Corollary 14, but the proof is different.)

**Corollary 19.** *Let $s$ be a real number, $m > s - \frac{1}{2}$ a nonnegative integer, and $x_0$ any point of $\Omega$. Then a holomorphic function $f$ belongs to $W^s_{\mathrm{hol}}(\Omega)$ if and only if any of the following quantities is finite, and the square root of each of these quantities gives an equivalent norm in $W^s_{\mathrm{hol}}(\Omega)$:*

(a) $\sum_{|\nu| \le m} \frac{|\nu|!}{\nu!} \|\partial^\nu f\|^2_{L^2(\Omega, |r|^{2m-2s})}$;

(b) $\sum_{|\nu| = m} \frac{m!}{\nu!} \|\partial^\nu f\|^2_{L^2(\Omega, |r|^{2m-2s})} + \|f\|^2_{L^2(\Omega, |r|^{2m-2s})}$;

(c) $\sum_{|\nu| = m} \frac{m!}{\nu!} \|\partial^\nu f\|^2_{L^2(\Omega, |r|^{2m-2s})} + \sum_{|\nu| < m} |\partial^\nu f(x_0)|^2$;

(d) $\sum_{j=0}^m \|\mathcal{D}^j f\|^2_{L^2(\Omega, |r|^{2m-2s})}$ *if* $m > 2s - 1$;

(e) $\|\mathcal{D}^m f\|^2_{L^2(\Omega, |r|^{2m-2s})} + \|f\|^2_{L^2(\Omega, |r|^{2m-2s})}$ *if* $m > 2s - 1$.

Note that (d) is precisely the result of Beatrous mentioned in the Introduction; our proof is totally different from the one in [1].

*Proof.* (a) Setting as usual $f = \mathbf{K}u$, we have seen in the proof of Theorem 8 that the expression in (a) equals

$$\sum_{|\nu| \le m} \binom{|\nu|}{\nu} \| |r|^{m-s} \mathbf{K} Z^\nu u \|^2_{L^2(\Omega)} = \langle Tu, u \rangle_{H^2(\partial\Omega)}$$

where

$$T = \mathbf{\Pi} \sum_{|\nu| \le m} \binom{|\nu|}{\nu} Z^{*\nu} \mathbf{K}^* |r|^{2m-2s} \mathbf{K} Z^\nu \mathbf{\Pi}$$

$$= \mathbf{\Pi} \sum_{|\nu| \le m} \binom{|\nu|}{\nu} Z^{*\nu} \Lambda_{|r|^{2m-2s}} Z^\nu \mathbf{\Pi} = T_Q$$

with

$$Q = \sum_{|\nu| \le m} \binom{|\nu|}{\nu} Z^{*\nu} \Lambda_{|r|^{2m-2s}} Z^\nu$$

a $\Psi$DO of order $2m + 2s - 2m - 1 = 2s - 1$ and with symbol which satisfies

$$\sigma(Q)|_\Sigma = \sum_{|\nu| = m} \binom{m}{\nu} \Big| \Big(\frac{\|\xi\|}{\|\partial r\|} \partial r\Big) \Big|^{2\nu} \frac{\Gamma(2m-2s+1)}{\|\xi\|^{2m-2s+1}} \|\partial r\|^{2m-2s}$$

$$= \Gamma(2m - 2s + 1) \|\xi\|^{2s-1} \|\partial r\|^{2m-2s} > 0.$$

Since the operator $T$ is nonnegative self-adjoint (being a sum of expressions of the form $Q^*Q$, with $Q = |r|^{m-s} \mathbf{K} Z^\nu : H^2(\partial\Omega) \to L^2(\Omega)$ densely-defined and closed) and positive (since the term $|\nu| = 0$ vanishes only when $u \equiv 0$), the assertion therefore follows by part (a) of Theorem 17.

(b) and (c) are proved in the same way.

(d) Setting again $u = \gamma f$, the proof of Theorem 9 shows that the expression in (d) equals

(58)
$$\sum_{j=0}^m (-1)^j \langle \mathbf{K}^* \mathcal{D}^{*j} |r|^{2m-2s} \mathcal{D}^j \mathbf{K} u, u \rangle_{\partial\Omega}.$$

Here $\mathcal{D}$ is understood as the operator with maximal domain (i.e. $f \in \text{dom}\,\mathcal{D}$ if $f \in L^2_{\text{hol}}(\Omega)$ and $\mathcal{D}f$, taken in the distributional sense, belongs to $L^2(\Omega)$), and $-\mathcal{D}^*$ is its Hilbert-space adjoint; in view of the hypothesis $m > 2s - 1$, $|r|^{2m-2s}\mathcal{D}^j\mathbf{K}u$ vanishes to order $m$, hence belongs to $\text{dom}\,\mathcal{D}^{*j}$. Thus

$$T := \sum_{j=0}^{m} (-1)^j \mathbf{K}^* \mathcal{D}^{*j} |r|^{2m-2s} \mathcal{D}^j \mathbf{K}$$

is a positive (since the term $j = 0$ in (58) vanishes only if $u \equiv 0$) self-adjoint operator on $L^2(\partial\Omega)$. Further, we have also seen in the proof of Theorem 9 that

$$\mathcal{D}^{*j} |r|^{2m-2s} \mathcal{D}^j f = \sum_{|\nu| \leq j} \phi_{j\nu} |r|^{2m-j-2s} \partial^\nu f$$

with $\phi_{j\nu} \in C^\infty(\overline{\Omega})$, and $\phi_{j\nu}|_{\partial\Omega} = (-1)^j \binom{j}{\nu} (\overline{\partial r})^\nu \frac{\Gamma(2m-2s+1)}{\Gamma(2m-j-2s+1)} \|\partial r\|^{2j}$ if $|\nu| = j$. It follows that

$$\mathbf{K}^* \mathcal{D}^{*j} |r|^{2m-2s} \mathcal{D}^j \mathbf{K} = \sum_{|\nu| \leq j} \Lambda_{\phi_{j\nu} |r|^{2m-j-2s}} Z^\nu$$

is a $\Psi$DO of degree $2s + 2j - 2m - 1$. Thus $\Pi T|_{H^2(\partial\Omega)}$ is a generalized Toeplitz operator of order $2s - 1$, with symbol

$$\frac{\Gamma(2m-2s+1)}{\Gamma(m-2s+1)} \|\partial r\|^{2m} \sum_{|\nu| = m} \binom{m}{\nu} \frac{\Gamma(m-2s+1)}{\|\xi\|^{m-2s+1}} \|\partial r\|^{m-2s} (\overline{\partial r})^\nu \left( \frac{\|\xi\|}{\|\partial r\|} \partial r \right)^\nu$$
$$= \Gamma(2m-2s+1) \|\partial r\|^{4m-2s} \|\xi\|^{2s-1} > 0,$$

and the claim follows by part (a) of Theorem 17.

Finally, (e) is proved in the same way as (d). $\quad\square$

*Remark 20.* Incidentally, in (a) and (d) in the last corollary, part (b) of Theorem 17 also gives another proof of Theorems 8 and 9. However, it seemed more transparent to give the direct proofs in Section 5. $\quad\square$

To some extent, the operator $\Lambda_0^{-2}$ on $\partial\Omega$ plays a similar role as the shifted Laplacian $I + \Delta$ on $\mathbf{R}^n$, and the expression $\langle \Lambda_0^{1-2s} u, u \rangle$ is an analogue of the Sobolev norm $\|\cdot\|'$ from (13). For this reason, the following assertion forms the last part of Theorem B which remains to be proved.

**Corollary 21.** *For $s \in \mathbf{R}$ let $\mathcal{H}_s$ be the completion of $C^\infty_{\text{hol}}(\overline{\Omega})$ with respect to the norm*

$$\|f\|'_s := \langle (T_{\Lambda_0})^{1-2s} \gamma f, \gamma f \rangle_{H^2(\partial\Omega)}^{1/2}.$$

*Then $\mathcal{H}_s = W^s_{\text{hol}}(\Omega)$, with equivalent norms, and the reproducing kernels of $\mathcal{H}_s$ are of the form (45) for $\alpha = -2s$, with the leading term given by (46) for $Q = \Lambda_0^{-2s}$, i.e.*

$$a(x,x) = \frac{\Gamma(n-2s+1)\, J[\rho](x)}{\pi^n \|2\partial r(x)\|^{2s}}$$

for $x \in \partial\Omega$, where $\Gamma(n - 2s + 1)$ is to be replaced by $(-1)^{k+1}/k!$ if $n - 2s + 1 = -k$ is a nonpositive integer.

*Proof.* Apply Theorem 17 to $T = (T_{\Lambda_0})^{1-2s}$, noting that $\Lambda_0 = \mathbf{K}^*\mathbf{K}$ is a positive self-adjoint $\Psi$DO of order $-1$ with symbol $\|2\xi\|^{-1}$. $\square$

Finally, we can use Proposition 16 to prove Theorem C.

*Proof of Theorem C.* As in the proof of Proposition 16, let $Q$ be a $\Psi$DO of order $-1$ such that $Q^* = Q$, $\sigma(Q) > 0$, $T_Q = T_{\Lambda_0}$ and $Q\Pi = \Pi Q$. Thus $\sigma(Q)|_\Sigma = \|2\xi\|^{-1}$, the operator $T$ from the preceding corollary is just $Q^{1-2s}$ restricted to $H^2(\partial\Omega)$, and the reproducing kernels $K^{(s)}(x, y) \equiv K_y^{(s)}(x)$ occurring there are just

$$K_y^{(s)} = Q^{2s-1} S_y$$

where $S_y(x) \equiv S(x, y)$ is the Szegő kernel (cf. Remark 18). Note that $S_y \in C_{\text{hol}}^\infty(\partial\Omega) \subset W_{\text{hol}}^{1-2\operatorname{Re}s}(\partial\Omega) = \operatorname{dom} Q^{2s-1}$ by the usual ellipticity argument (cf. the proof of Proposition 16).

Let $\lambda_j$ and $P_j$, $j = 1, 2, \ldots,$ be the eigenvalues and their spectral projections, respectively, for the (compact positive self-adjoint) operator $Q|_{H^2(\partial\Omega)}$. Then by the spectral theorem

$$K_y^{(s)} = \sum_j \lambda_j^{2s-1} P_j S_y,$$

and thus

$$\begin{aligned}
K^{(s)}(x, y) &= \langle K_y^{(s)}, S_x \rangle_{H^2(\partial\Omega)} \\
&= \sum_j \lambda_j^{2s-1} \langle P_j S_y, S_x \rangle \\
&= \sum_j \lambda_j^{2s-1} \langle P_j S_y, P_j S_x \rangle
\end{aligned}$$

(59)

(since $P_j$ is a projection). Now

$$\begin{aligned}
|K^{(s)}(x, y)| &\leq \sum_j \lambda_j^{2\operatorname{Re}s-1} \|P_j S_y\| \|P_j S_x\| \\
&\leq \Big( \sum_j \lambda_j^{2\operatorname{Re}s-1} \|P_j S_y\|^2 \Big)^{1/2} \Big( \sum_j \lambda_j^{2\operatorname{Re}s-1} \|P_j S_x\|^2 \Big)^{1/2} \\
&= \sqrt{K^{(\operatorname{Re}s)}(x, x) \, K^{(\operatorname{Re}s)}(y, y)} \\
&= \|K_x^{(\operatorname{Re}s)}\|_{\operatorname{Re}s} \|K_y^{(\operatorname{Re}s)}\|_{\operatorname{Re}s} \\
&= \|Q^{\operatorname{Re}s-1/2} S_x\|_{H^2(\partial\Omega)} \|Q^{\operatorname{Re}s-1/2} S_y\|_{H^2(\partial\Omega)} < \infty,
\end{aligned}$$

since $S_x \in \operatorname{dom} Q^s$ for all $s \in \mathbf{C}$ and $x \in \Omega$, as noted above. Thus the series (59) converges for any $s \in \mathbf{C}$ and $x, y \in \Omega$. Also, since $K^{(\operatorname{Re}s)}(x, x)$ — being the restriction to the diagonal of a sesqui-analytic function — is continuous (even real-analytic) on $\Omega$, it follows from the third line in the last chain of inequalities that the convergence is even uniform for $x, y$ in compact subsets of $\Omega$. Finally, since $0 < \lambda_j \leq \|Q\| < \infty$, the fact that (59) converges absolutely for some $s_0 \in \mathbf{R}$ implies that it also converges, and uniformly so, in the strip $s_0 \leq \operatorname{Re} s \leq s_0 + \kappa$ for any $\kappa > 0$.

In conclusion, the series (59) converges uniformly for $(x, y, s)$ in compact subsets of $\Omega \times \Omega \times \mathbf{C}$, and thus defines a holomorphic functions of $(x, \overline{y}, s)$ there, which coincides with the Sobolev-Bergman kernels $K^{(s)}(x, y)$ from Corollary 21 for $s \in \mathbf{R}$. This completes the proof of Theorem C.  □

We remark that, in fact, the usual argument involving the parametrix and the formulas (45) and (46) shows that the description of the boundary singularities of $K^{(s)}(x, y)$ given in Corollary 21 remains in force even for complex $s$.

## 8. Concluding remarks

**8.1 Equivalence.** It has been alluded to, at several places above, that the reproducing kernels of the same space with respect to two equivalent norms may have very different boundary singularities. Here is an example.

Take the unit disc $\mathbf{D} = \{z \in \mathbf{C} : |z| < 1\}$. With respect to the normalized Lebesgue measure, the monomials $z^k$, $k = 0, 1, 2, \ldots$, form an orthogonal basis, with

$$(60) \qquad \|z^k\|^2 = \frac{1}{k+1},$$

and the corresponding reproducing kernel is the traditional Bergman kernel

$$K(x, y) = \sum_{k=0}^{\infty} \frac{(x\overline{y})^k}{\|z^k\|^2} = \frac{1}{(1 - x\overline{y})^2}.$$

Introduce another scalar product in $L^2_{\mathrm{hol}}(\mathbf{D})$ by letting

$$(61) \qquad \langle z^j, z^k \rangle = \frac{\delta_{jk}}{k + 1 + a_k}$$

where

$$(62) \qquad a_k \geq 0, \qquad \sup_k \frac{a_k}{k+1} < \infty.$$

Then the corresponding norms are clearly equivalent. The reproducing kernel with respect to (61) is given by

$$K'(x, y) = \sum_{k=0}^{\infty} (k + 1 + a_k)\,(x\overline{y})^k.$$

Now choosing

$$a_k = \frac{\Gamma(k + \frac{3}{2})}{k!\,\Gamma(\frac{3}{2})} \simeq \sqrt{k+1}$$

the corresponding reproducing kernel is

$$K'(x, y) = (1 - x\overline{y})^{-2} + (1 - x\overline{y})^{-3/2},$$

which is not of the form (10). More generally,

$$a_k = \frac{\Gamma(k+\beta)}{k!\Gamma(\beta)} \simeq (k+1)^{\beta-1},$$

which satisfies (62) for any $0 < \beta \leq 2$, produces

$$K'(x,y) = (1 - x\overline{y})^{-2} + (1 - x\overline{y})^{-\beta}.$$

Taking

$$a_0 = 0, \qquad a_k = \log k \quad (k \geq 1),$$

we get

$$K'(x,y) = \frac{1}{(1 - x\overline{y})^2} + \sum_{k=1}^{\infty} (x\overline{y})^k \log k$$
$$= \frac{1}{(1 - x\overline{y})^2} + \frac{1}{1 - x\overline{y}} \log \frac{1}{1 - x\overline{y}} + O\left(\frac{1}{1 - x\overline{y}}\right).$$

It is not difficult to construct examples of even wilder boundary singularities.

## 8.2 Complex powers of weights.

The various kernels occurring for a given Sobolev-Bergman space $W_{\text{hol}}^s(\Omega)$ in Theorems A, 8, 9, etc., almost never coincide. In particular, the kernels from Theorem C are different from the ones from Theorem A (with $\alpha = -2s$, $s < \frac{1}{2}$), even though $L_{\text{hol}}^2(\Omega, |r|^{-2s}) = W_{\text{hol}}^s(\Omega)$ as spaces, with equivalent norms. Thus Theorem C does not imply that the kernels $K^{(\alpha)}(x,y)$ of $L_{\text{hol}}^2(\Omega, |r|^\alpha)$ from Theorem A can be holomorphically continued to $\alpha \in \mathbf{C}$, or at least to $\text{Re}\,\alpha > -1$.

CONJECTURE. $K^{(\alpha)}(x,y)$ extends to a holomorphic function of $x, \overline{y}, \alpha$ on $\Omega \times \Omega \times \{\text{Re}\,\alpha > -1\}$.

Of course, using Proposition 2, this is tantamount to having an analytic continuation of

$$(T_{|r|^\alpha})^{-1} K_x\,(y), \qquad \text{Re}\,\alpha > -1.$$

While there is no problem with the holomorphy of $\alpha \mapsto T_{|r|^\alpha}$ in an appropriate sense (cf. [21], Chapter VII), the difficulty lies with taking the inverse: for $\alpha \notin \mathbf{R}$, there seems to be no reason to expect $T_{|r|^\alpha}$ to be injective, or $K_x$ to be in its range. Another possible line of attack — though yielding an analytic continuation only to a small neighbourhood of the positive real half-axis — would be to estimate the derivatives of the functions $\alpha \mapsto T_{|r|^\alpha}^{-1}$ on $\alpha > -1$, that is, expressions of the form

$$\langle T_{|r|^\alpha}^{-1} T_{|r|^\alpha (\log|r|)^{k_1}} T_{|r|^\alpha}^{-1} T_{|r|^\alpha (\log|r|)^{k_2}} T_{|r|^\alpha}^{-1} \cdots T_{|r|^\alpha (\log|r|)^{k_m}} T_{|r|^\alpha}^{-1} K_y, K_x \rangle,$$

but this does not seem to be any easier.

It seems even likely that $K^{(\alpha)}(x,y)$ extends to a holomorphic function of $x, \overline{y}, \alpha$ on all of $\Omega \times \Omega \times (\mathbf{C} \setminus \mathbf{Z}_{<0})$, with simple poles at the negative integers.

## 8.3 Boundary invariants.

The analysis of boundary singularities of various reproducing kernels can be used towards obtaining important $CR$-invariants of the domain, cf. e.g. [18], [15], and the references there. (This was, in fact, our original motivation for undertaking this study.) It would be of interest to know if suitable versions of the various Sobolev-Bergman kernels above can be useful in this regard.

**8.4 Logarithmic weights.** Especially from the point of view of the last remark, it would be desirable to extend the results of this paper to generalized Toeplitz operators with "logarithmic" terms in the symbol, i.e. to reproducing kernels of spaces like $L^2_{\text{hol}}(\Omega, |r|^\alpha |\log|r||^\beta)$, $\alpha > -1$, $\beta \in \mathbf{R}$. The reason is that the various kernels above — whose construction relied on quantities like the defining function $r$ — are manifestly *not* invariant under biholomorphisms: there does not exist any holomorphic-invariant recipe for a defining function, for instance (see Hirachi and Komatsu [17], §5.1), and similarly the definitions of the operators $\mathcal{D}$ and $\overline{\mathcal{D}}$ use the ambient Euclidean structure rather than the intrinsic geometry of $\Omega$. A way out of this might be replacing all those objects by suitable holomorphically-invariant ones, e.g. the defining function $-r$ by the solution $u$ of the Monge-Ampère equation $J[u] = 1$, or the term $e^g$ in (7) by the Bergman invariant $\beta_K := K^{-1} \det[\partial\overline{\partial} \log K]$ ($K$ being the ordinary Bergman kernel), or $\sum_{|\alpha|=m} |\partial^\alpha f|^2$ by $\widetilde{\Delta}^m |f|^2$ (for $f$ holomorphic), where $\widetilde{\Delta}$ is the Laplace-Beltrami operator with respect to some invariant metric (Bergman, Poincaré, etc.). (It is not completely clear what should be the "invariant" substitute for $\mathcal{D}$.) All the objects just mentioned — $u$, $\beta$, as well as $\widetilde{\Delta}$ — have a logarithmic singularity at $\partial\Omega$ of some sort (see [24], [11], [23]). For ΨDOs, some results on these "logarithmic symbols" do exist (see e.g. [32]); however, their analogues for generalized Toeplitz operators seem to be currently out of reach.

**8.5 Hermite FIOs.** Throughout this paper, we have largely followed the exposition of generalized Toeplitz operators in the paper [6] and in the appendix of the book [7], which rely on FIOs with complex-valued phase functions. The rest of [7] uses the more general "FIOs of Hermite type" instead, which seem to have the drawback that, apparently, for them one really needs to assume that the order of $T_Q$ is an integer or a half-integer. (What breaks down is the "parity" argument on pp. 75–76 in [7], needed to show that certain symbols are classical, instead of having an asymptotic expansion where the degrees of homogeneity go down by steps of only $\frac{1}{2}$ instead of 1; this is in turn needed in the proof of the equivalences (31) and, hence, of the properties (P1)–(P7) (cf. Proposition 2.3 in [7]). For FIOs with complex-valued phase functions instead of Hermite-type FIOs, things work fine. (The reason is that degrees are not counted in the same manner for Hermite operators and for FIOs: e.g. if $Q$ is of degree 0 and its symbol vanishes to order 1 on $\Sigma$, then $\Pi Q \Pi$ is of degree $-1/2$ as a Hermite operator, but usually still of degree 0 as a FIO. The author is grateful to Louis Boutet de Monvel for this clarification.))

MATHEMATICS INSTITUTE, SILESIAN UNIVERSITY AT OPAVA, NA RYBNÍČKU 1, 74601 OPAVA, CZECH REPUBLIC and MATHEMATICS INSTITUTE, ŽITNÁ 25, 11567 PRAGUE 1, CZECH REPUBLIC
*E-mail address:* englis@math.cas.cz